\documentclass[11pt,reqno,oneside]{article}

\usepackage{multirow}
\usepackage{cancel}
\usepackage{pstricks,pst-node,pst-tree}
\usepackage{amsmath}
\usepackage{amssymb}
\usepackage[lmargin=.55in, rmargin=.55in, tmargin=1in, bmargin=1in]{geometry}
\usepackage{xcolor}
\usepackage{graphicx}
\usepackage[font=small]{caption}
\usepackage[labelformat = empty,position=top]{subcaption}
\usepackage{pstricks}
\newtheorem{theo}{Theorem}
\newtheorem{lemm}[theo]{Lemma}
\newtheorem{prop}[theo]{Proposition}

\definecolor{myblue}{rgb}{0,0.1,0.9}
\definecolor{myred}{rgb}{0.9,0.1,0}

\renewcommand{\geq}{\geqslant}
\renewcommand{\leq}{\leqslant}

\title{{\bf The distributions under two species-tree models \\ of the total number of ancestral configurations \\ for matching gene trees and species trees}}

\author{Filippo Disanto\thanks{Department of Mathematics, University of Pisa, Pisa 56126, Italy. Email: filippo.disanto@unipi.it.},
Michael Fuchs\thanks{Department of Mathematical Sciences, National Chengchi University, Taipei 116, Taiwan. Email: mfuchs@nctu.edu.tw.},
Chun-Yen Huang\thanks{Department of Applied Mathematics, National Chiao Tung University, Hsinchu 300, Taiwan.},
Ariel R.~Paningbatan\thanks{Institute of Mathematics, University of the Philippines Diliman, Quezon City 1101, Philippines. Email: arpaningbatan@math.upd.edu.ph.},
Noah A.~Rosenberg\thanks{Department of Biology, Stanford University, Stanford, CA 94305, USA. Email: noahr@stanford.edu.}}

\begin{document}

\maketitle

\begin{abstract}
Given a gene-tree labeled topology $G$ and a species tree $S$, the \emph{ancestral configurations} at an internal node $k$ of $S$ represent the combinatorially different sets of gene lineages that can be present at $k$ when all possible realizations of $G$ in $S$ are considered. Ancestral configurations have been introduced as a data structure for evaluating the conditional probability of a gene-tree labeled topology given a species tree, and their enumeration assists in describing the complexity of this computation. In the case that the gene-tree labeled topology $G=t$ matches that of the species tree $S$, by techniques of analytic combinatorics, we study distributional properties of the \emph{total} number of ancestral configurations measured across the different nodes of a random labeled topology $t$ selected under the uniform and the Yule probability models. Under both of these probabilistic scenarios, we show that the total number $T_n$ of ancestral configurations of a random labeled topology of $n$ taxa asymptotically follows a lognormal distribution. Over uniformly distributed labeled topologies, the asymptotic growth of the mean and the variance of $T_n$ are found to satisfy $\mathbb{E}_{\rm U}[T_n] \sim 2.449 \cdot 1.333^n$ and $\mathbb{V}_{\rm U}[T_n] \sim 5.050 \cdot 1.822^n$, respectively. Under the Yule model, which assigns higher probabilities to more balanced labeled topologies, we obtain the mean $\mathbb{E}_{\rm Y}[T_n] \sim 1.425^n$ and the variance $\mathbb{V}_{\rm Y}[T_n] \sim 2.045^n$.
\end{abstract}

\begin{itemize}
\item[] {\bf Keywords:} analytic combinatorics, ancestral configurations, gene trees, phylogenetics, species trees.
\item[] {\bf Mathematics subject classification (2010):} 05A15 $\cdot$ 05A16 $\cdot$ 05C05 $\cdot$ 92D15
\end{itemize}

\section{Introduction}

Ancestral configurations are lists that describe for a given gene-tree topology $G$ and a species tree $S$ the sets of gene lineages that can be present at a given node of $S$ (Fig.~\ref{configa}). They have been introduced by Wu \cite{Wu12} as a data structure in the calculation of the probabilities of gene-tree topologies conditional on species trees under the multispecies coalescent model. In particular, for a given species tree $S$, Wu's algorithm ``STELLS'' evaluates the probability of a gene-tree topology $G$ by recursively computing the probabilities of the ancestral configurations of $G$ at the different nodes of $S$, proceeding from the tips towards the root of $S$~\cite{Wu12}. The running time of STELLS depends on the total number $c(G, S)$ of ancestral configurations of $G$ in $S$, that is, on the sum of the number of ancestral configurations of $G$ across the different nodes of $S$.

If the topology $G=t$ of the gene tree matches the topology of species tree $S$, then the total number of configurations of $G$ in $S$ becomes a function $c(G, S) = c(t)$ of $t$, whose behavior over tree families of increasing size can be analyzed by tools of enumerative and analytic combinatorics~\cite{FlajoletAndSedgewick09}. In initial studies~\cite{DisantoEtAl20, DisantoAndRosenberg17}, by examining the number of ancestral configurations at the root of a randomly selected topology $t$ with number of leaves $n$, we derived theorems that determine the exponential growth---indicated here by the symbol ``$\bowtie$''---of the mean $\mathbb{E}[c]$ and the variance $\mathbb{V}[c]$ of the total number of configurations. In particular, we found that the mean grows exponentially like $\mathbb{E}_{\rm U}[c] \bowtie 1.333^n$ and $\mathbb{E}_{\rm Y}[c] \bowtie 1.425^n$ for uniformly and Yule-distributed topologies of size $n$, respectively. Under the same distributions, the exponential growth of the variance satisfies $\mathbb{V}_{\rm U}[c] \bowtie 1.822^n $ and $\mathbb{V}_{\rm Y}[c] \bowtie 2.045^n$.


These results, however, do not fully characterize the \emph{sub}exponential growth of the mean and variance of the random variable $c$, and the problem of describing the asymptotic distribution of the total number of configurations has remained open. Here, we solve these problems by using generating functions to count the total number of ancestral configurations in random tree topologies. Surprisingly, we find that up to a constant factor---which we calculate exactly---the exponential growth of $\mathbb{E}_{\rm U}[c]$, $\mathbb{E}_{\rm Y}[c]$, $\mathbb{V}_{\rm U}[c]$, and $\mathbb{V}_{\rm Y}[c]$ determines the full asymptotic behavior of the associated quantities.
More precisely, for random topologies of increasing size $n$ selected under the uniform and Yule distributions, we show that $\mathbb{E}_{\rm U}[c] \sim 2.449 \cdot 1.333^n$, $\mathbb{E}_{\rm Y}[c] \sim 1.425^n$, $\mathbb{V}_{\rm U}[c] \sim 5.050 \cdot 1.822^n$, and $\mathbb{V}_{\rm Y}[c] \sim 2.045^n$.
Furthermore, under both the uniform and Yule models, we obtain an asymptotic lognormal distribution of the total number of ancestral configurations. We study the correlation between the total number of configurations and the closely related number of root ancestral configurations.

Our approach uses standard techniques of analytic combinatorics for deriving the asymptotic growth of integer sequences coupled with a key observation that enables the study of distributional properties of the number of ancestral configurations of random tree topologies selected under the uniform and Yule models by equivalently considering uniformly distributed classes of plane trees often known as \emph{Catalan trees} and \emph{increasing binary trees}. The results contribute to the enumerative study of combinatorial structures in the relationship between gene trees and species trees~\cite{AlimpievandRosenberg21, DegnanAndSalter05, DisantoAndMunarini19, DisantoAndRosenberg15, DisantoAndRosenberg16, DisantoAndRosenberg19JMB, DisantoAndRosenberg19:bmb, HimwichAndRosenberg20, Rosenberg07:jcb, Rosenberg13:tcbb, Rosenberg19, RosenbergAndDegnan10}, and they can assist in relating the complexity of algorithms for computing gene-tree probabilities with ancestral configurations to algorithms that use an evaluation based on other structures~\cite{DegnanAndSalter05, TruszkowskiEtAl21, Wu16}.


\section{Preliminaries}

We start with some definitions, preliminary results, and basic principles of enumerative combinatorics. In Section \ref{labetopolo}, we introduce labeled topologies and their uniform and Yule probability distributions. In Section \ref{analyticomb}, we present generating function techniques for use in analyzing the asymptotic growth of integer sequences.

\subsection{Labeled and unlabeled topologies}
\label{labetopolo}

\begin{figure}
\begin{center}
\includegraphics*[scale=0.74,trim=0 0 0 0]{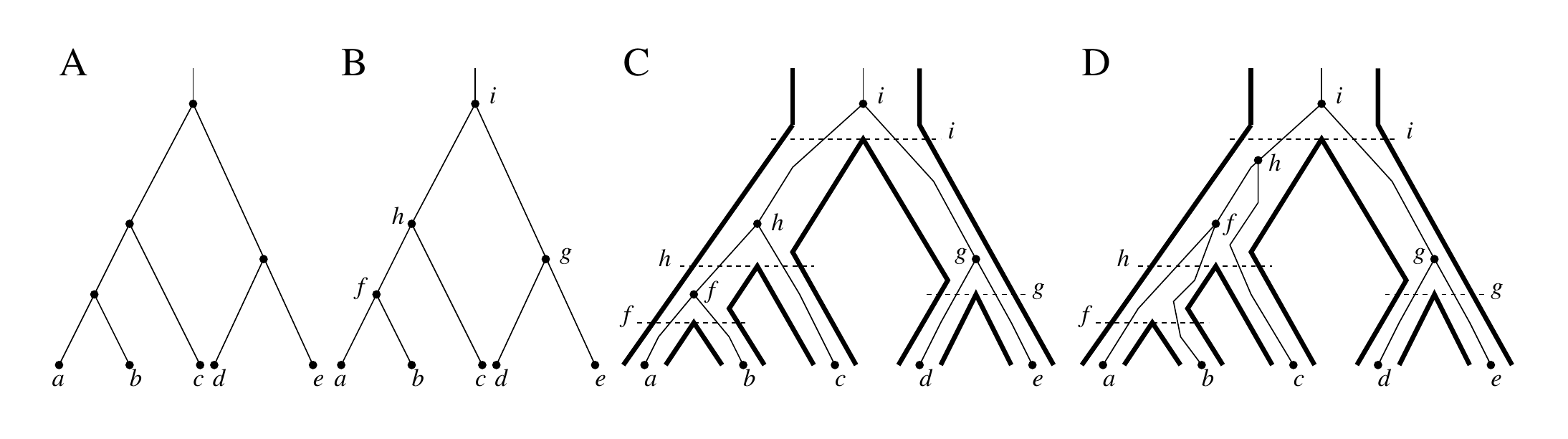}
\end{center}
\vspace{-.7cm}
\caption{{\small Labeled topologies, gene trees, and species trees. {\bf (A)} A labeled topology of size 5. {\bf(B)} The labeled topology in (A) with its internal nodes labeled. We identify each edge of the tree by its immediate descendant node; for example, lineage $h$ results from the coalescence of lineages $c$ and $f$. {\bf (C)} A possible realization (thin lines) of the gene-tree labeled topology of (A) in a species tree with a matching labeled topology (thick lines). The ancestral configuration at species-tree node $i$ is $\{ g, h \}$. The configuration at node $h$ is $\{ c, f \}$. {\bf (D)} A different realization of the gene-tree labeled topology in (A) in a matching species tree. The ancestral configurations at species-tree nodes $i$ and $h$ are $\{ g, h \}$ and $\{ a, b, c \}$, respectively.
}} \label{configa}
\end{figure}

A \emph{labeled topology} $t$, also called a \emph{phylogenetic tree}, of size $|t|=n$ is a binary rooted tree whose $n$ external nodes---its leaves---possess distinct labels, often for small $n$ the first $n$ letters alphabetically (Fig.~\ref{configa}A). Labeled topologies are \emph{non-plane}, or \emph{unordered}, in the sense that each pair of child nodes carries no left-right orientation; we obtain the same labeled topology by transposing the two subtrees stemming from an internal node.

It is convenient to denote the internal nodes of a labeled topology $t$ by letters different from those associated with the leaves (Fig.~\ref{configa}B). We identify each edge of $t$ by (the label of) its immediate descendant node, i.e., by the node closer to the leaves that is adjacent to the edge. We describe the descendant--ancestor order relation defined over the set of nodes of $t$ by the symbol $\preceq$. More precisely, for distinct nodes $x$ and $y$, we write $x \prec y$ in $t$ if $y$ is a node belonging to the path connecting node $x$ to the root of $t$. The subtree of $t$ rooted at node $k$, which contains those nodes $x$ of $t$ with $x \preceq k$, is denoted by $t^k$. Hence, in particular, $t^k=t$ if $k$ is the root of $t$, and $t^k= \bullet_{k}$ if $k$ is a leaf of $t$, where $\bullet_{k}$ is a subtree that contains only node $k$.

By removing labels from a labeled topology $t$, we obtain the tree \emph{shape} or \emph{unlabeled topology} underlying $t$. Unlabeled topologies are also called \emph{Otter trees}~\cite{Otter48}; for increasing numbers of leaves $n \geq 1$, they are enumerated by the Wedderburn-Etherington numbers, $1, 1, 1, 2, 3, 6, 11, 23, 46, \ldots$~\cite{Felsenstein04, Wedderburn22}.

Distinct labeled topologies $t_1$ and $t_2$ can possibly share the same unlabeled topology. For instance, in \emph{Newick format}, labeled topologies $t_1=((a,b),c)$ and $t_2=(a,(b,c))$ share unlabeled topology $((\bullet,\bullet),\bullet) = (\bullet,(\bullet,\bullet))$. The number ${\rm lab}(t)$ of labeled topologies with shape $t$ is obtained recursively by eq.~(22) of \cite{DisantoEtAl20},
\begin{equation}\label{lteq}
{\rm lab}(t) = {\rm lab}(t_L) \, {\rm lab}(t_R) \, {{|t|}\choose{|t_L|}}\frac{1}{1 + \delta_{t_L = t_R}}.
\end{equation}
Here, $t_L$ and $t_R$ are the two subtrees stemming from the root of $t$ (\emph{root subtrees}, for short) and $\delta_{t_L = t_R}$ is the Kronecker delta that equals 1 if $t_L$ and $t_R$ are the same unlabeled topology. We set ${\rm lab(t)} = 1$ if $|t|=1$.

Let $L_n$ denote the set of labeled topologies of size $n$. For $n \geq 2$, the cardinality of $L_n$ is $|L_n|=(2n-3)!! = 1 \times 3 \times 5 \times \ldots \times (2n-3)$~\cite{Felsenstein78, Felsenstein04}, which can be written
\begin{equation}\label{carciofo}
|L_n| =
\frac{(2n)!}{2^n(2n-1)n!}.
\end{equation}

Different probability models can be considered over the set $L_n$ of labeled topologies of fixed size $n$~\cite{Aldous01}. Under the uniform model, each labeled topology $t \in L_n$ has equal probability
$$\mathbb{P}_U[t] = \frac{1}{|L_n|} = \frac{2^{n} (2n-1) n!}{(2n)!}.$$
The Yule model is a generative model in which each lineage is equally likely to be the next to bifurcate forward in time, or equivalently, each pair of lineages is equally likely to be next to merge back in time~\cite{Harding71, Steel16, Yule25}. Many of its combinatorial features have been studied~\cite{ChangAndFuchs10, DisantoEtAl13, DisantoAndWiehe13, McKenzieAndSteel00, Rosenberg06:anncomb}; a labeled topology $t$ of $n$ leaves has probability
\begin{equation}\label{pyule}
\mathbb{P}_Y[t] = \frac{2^{n-1}}{n! \prod_{r=3}^n (r-1)^{d_r(t)} }
\end{equation}
under the Yule model, where $d_r(t)$ is the number of nodes of $t$ with $r$ descending leaves~\cite{Brown94, RosenbergAndKing23, SteelAndMcKenzie01}.

Owing to the product appearing in the denominator, under the Yule distribution, more balanced labeled topologies tend to have larger probabilities \cite{Harding71}. For example, among the labeled topologies of size $5$, the one depicted in Fig.~\ref{configa}A has maximal Yule probability $\mathbb{P}[ (((a,b),c),(d,e)) ] = \frac{1}{60}$; taking $((((a,b),c),d),e)$ and $(((a,b),(c,d)),e)$ as representative labeled topologies for their unlabeled shapes, we have $\mathbb{P}[ ((((a,b),c),d),e) ] = \frac{1}{180}$ and $\mathbb{P}[ (((a,b),(c,d)),e) ] = \frac{1}{90}$.

\subsection{Asymptotic growth and generating functions}
\label{analyticomb}

This article studies the growth of non-negative integer sequences. We use the following notation. For two sequences $(a_n)_{n \geq 0}$ and $(b_n)_{n \geq 0}$, we write $a_n \sim b_n$ when the ratio $b_n/a_n$ converges to $1$ for $n \rightarrow \infty$. If $a_n \sim b_n$, then we say that, asymptotically, sequences $a_n$ and $b_n$ have the same growth. The sequence $a_n$ is said to have exponential growth $k^n$ or, equivalently, to be of exponential order $k$, when $a_n \sim k^n s(n),$ where $k$ is a constant and $s(n)$ is a subexponential factor. We write $a_n \bowtie b_n$ if $a_n$ and $b_n$ have the same exponential growth.

The generating function of a sequence $(a_n)_{n \geq 0}$ is the power series $A(z) = \sum_{n=0}^{\infty} a_n z^n$. Multiplying $A(z)$ by a generating function $B(z) = \sum_{n=0}^{\infty} b_n z^n$ gives the generating function $A(z) \, B(z) = \sum_{n=0}^{\infty} \sum_{j=0}^{n} a_j z^j \cdot b_{n-j} z^{n-j} = \sum_{n=0}^{\infty} (\sum_{j=0}^{n} a_j b_{n-j}) z^n$, whose $n$th coefficient is the convolution $\sum_{j=0}^{n} a_j b_{n-j}$. Also, if $k$ is a constant, then $A(z) + B(z) = \sum_{n=0}^{\infty} (a_n + b_n) z^n$ and $k A(z) = \sum_{n=0}^{\infty} (k a_n) z^n$.

If $A(z) = \sum_{n=0}^{\infty} a_n z^n$ is considered as a function of the complex variable $z$, then the analysis of $A(z)$ near its singularities---the points in the complex plane where $A(z)$ is not analytic---can assist in the study of the asymptotic growth of the coefficients $a_n = [z^n] A(z)$. The simplest scenario is when $A(z)$ has a unique dominant singularity $\alpha > 0$, that is, when $z = \alpha > 0$ is the only singularity of $A(z)$ of smallest modulus. In this case, under the fairly general conditions of Theorems IV.7 and VI.4 of \cite{FlajoletAndSedgewick09}, the singular expansion $A_{\alpha}(z)$ of the generating function $A(z)$ at $z = \alpha$ determines the asymptotic growth of the coefficients $a_n$ as
\begin{equation}\label{prima}
a_n \sim [z^n] A_{\alpha}(z) \bowtie \alpha^{-n}.
\end{equation}
In other words, the $n$th coefficient of $A(z)$ has for increasing $n$ the same growth as the $n$th coefficient of the expansion $A_{\alpha}(z)$, where, in particular, $\frac{1}{\alpha}$ is the exponential order of sequence $a_n$. For instance, as given by Example II.19 of \cite{FlajoletAndSedgewick09}, $L(z)= 1-\sqrt{1-2z}$ is the generating function associated with the sequence $|L_n|/n!$---where $L_n$ is the number of labeled topologies of size $n$ (eq.~\ref{carciofo}). The dominant singularity of $L(z)$ is $\alpha=\frac{1}{2}$, and indeed, in agreement with eq.~(\ref{prima}), we have $|L_n|/n! = {{2n}\choose{n}}/[2^n (2n-1)] \bowtie 2^n = \alpha^{-n}$.

In the following sections, we apply the asymptotic relation in eq.~(\ref{prima}) to generating functions $A(z)$ with a unique dominant singularity $\alpha > 0$ and singular expansion given by either
$A_{\alpha}(z) = 1/({1-\frac{z}{\alpha}})$
or
$A_{\alpha}(z) = k_1 - k_2 \sqrt{1 - {z}/{\alpha}},$
where $k_1$ and $k_2 > 0$ are constants. Using the equivalence in eq.~(\ref{prima}) with $ [z^n] \big(1/(1 - z)\big) = 1 $ and $[z^n] \big(- \sqrt{1 - z} \big) \sim 1/(2 \sqrt{\pi n^3})$~(see p.~388 of \cite{FlajoletAndSedgewick09}), these expansions yield
\begin{eqnarray}
\label{pillo1}
A_{\alpha}(z) & = & \frac{1}{1-\frac{z}{\alpha}} \quad \Rightarrow \quad a_n \sim \alpha^{-n} \\
\label{pillo2}
A_{\alpha}(z) & = & k_1 - k_2 \sqrt{1 - \frac{z}{\alpha}} \quad \Rightarrow \quad a_n \sim k_2 \, \frac{\alpha^{-n}}{2 \sqrt{\pi n^3}}.
\end{eqnarray}

If $A'(z) = \sum_{n=0}^{\infty} n a_n z^{n-1}$ and $\int_0^z A(t) \, dt = \sum_{n=0}^{\infty} [{a_n}/(n+1)] z^{n+1}$ are generating functions obtained by differentiating and integrating generating function $A(z)$, then we use Theorems VI.8 and VI.9(ii) of \cite{FlajoletAndSedgewick09} to calculate the singular expansions of $A'(z)$ and $\int_0^z A(t) \, dt$ at their dominant singularity $\alpha$ as
\begin{eqnarray}\label{diffe}
A'(z) & \stackrel{z \rightarrow \alpha}{\sim} & \Big(A_{\alpha}(z)\Big)' \\
\label{inte}
\int_0^z A(t) \, dt & \stackrel{z \rightarrow \alpha}{\sim} & \int_0^z A_{\alpha}(t) \, dt.
\end{eqnarray}
In particular, both $A'(z)$ and $\int_0^z A(t) \, dt$ have the same dominant singularity $\alpha$ as $A(z)$, and their singular expansions are obtained by respectively differentiating and integrating the singular expansion $A_{\alpha}(z)$ of $A(z)$. Note that we apply formula eq.~(\ref{inte}) to generating functions $A(z)$ with singular expansion $ A_{\alpha}(z) = {1}/{(1-\frac{z}{\alpha})^2}$, in agreement with the hypothesis of case (i) of Theorem VI.9 of \cite{FlajoletAndSedgewick09}.

\section{Ancestral configurations for matching gene trees and species trees}

In this section, we define ancestral configurations for matching gene trees and species trees (Section \ref{dex}) and explain how distributional properties of the number of ancestral configurations can be equivalently studied over labeled topologies and over other tree families (Section \ref{models}). Section \ref{distrib} reviews results of \cite{DisantoEtAl20} and \cite{DisantoAndRosenberg17} on ancestral configurations at the root of randomly selected labeled topologies.

\subsection{Definitions and examples}\label{dex}
We first introduce gene trees, species trees, and realizations of a gene tree in a species tree. Following \cite{Wu12}, we define ancestral configurations for pairs of gene trees and species trees that share the same labeled topology.

\subsubsection{Gene trees and species trees}

A \emph{species-tree} labeled topology represents the evolutionary relationships of a set of populations or species identified with the leaves of the tree. A \emph{gene-tree} labeled topology describes the evolution of (genomic regions of) individuals sampled from a set of populations or species.

If individuals are sampled from the populations considered at the leaves of the species tree, then the gene tree can be viewed from a biological perspective as a set of gene lineages (Fig.~\ref{configa}C, thin lines) that have spread by evolutionary forces within the branching structure of the species tree (Fig.~\ref{configa}C, thick lines). We assume that exactly one individual is sampled for each population at the leaves of the species tree. We only examine pairs of \emph{matching} species trees and gene trees, that is, pairs with the same labeled topology (Fig.~\ref{configa}C and D).

For a fixed species tree, the same gene tree can result from different instances, or realizations, of the evolutionary process. In panels C and D of Fig.~\ref{configa}, the gene-tree labeled topology of panel A---with internal nodes denoted as in panel B---is depicted within the species tree as an outcome of two realizations that differ in the choice of the branches (edges) of the species tree where the coalescent events (internal nodes) of the gene tree occur. In particular, switching from panel C to panel D, we find gene-tree coalescent event $f$ in two different species-tree branches. In mathematical terms, a \emph{realization} of gene tree $G$ in species tree $S$ with matching labeled topology $S=G=t$ is a function $R$ mapping the set of internal nodes of $t$ onto itself such that two conditions hold: (i) for all internal nodes $k$, $k \preceq R(k)$, and (ii) for all internal nodes $k_1$ and $k_2$, $k_1 \preceq k_2 \Rightarrow R(k_1) \preceq R(k_2)$.
By identifying each species-tree branch by its immediate descendant node, the coalescent event corresponding to internal node $k$ of the gene tree $G$ is specified by the realization $R$ to occur in branch $R(k)$ of the species tree. For example, the realization that encodes the evolutionary scenario in Fig.~\ref{configa}C is $R(k) = k$ for all $k \in \{ f, g, h, i \}$, whereas in Fig.~\ref{configa}D, the realization instead has $R(f) = R(h) = h$, $R(g) = g$, and $R(i) = i$.
\subsubsection{Ancestral configurations}\label{defex}

When species trees are equipped with branch lengths that measure the time separating pairs of adjacent nodes, the conditional probability of a gene-tree labeled topology for a given species tree can be calculated under the multispecies coalescent model~\cite{DegnanAndSalter05}. Ancestral configurations were introduced by Wu \cite{Wu12} as a data structure for the recursive calculation of this conditional probability, with each node of the species tree associated with a given set of ancestral configurations depending on the possible realizations of the gene tree. At each step, Wu's algorithm computes the probability under the coalescent model that an ancestral configuration at a given node of the species tree has ``evolved'' from the ancestral configurations at its child nodes, proceeding recursively from the leaves to the root. The cost of Wu's algorithm is affected by the total number of ancestral configurations measured across all nodes of the species tree.

In our setting, where the gene-tree labeled topology matches that of the species tree, ancestral configurations are defined as follows. Given a gene-tree labeled topology $G=t$ and a matching species tree $S$, let $R$ be a realization of $G$ in $S$.
For a given node $k$ of $S$, consider the set $C(k) = C(k,R)$ of gene lineages (edges of $G$) present in $S$ at the time point right before node $k$, when time flows from the leaves toward the species-tree root.
The set $C(k)$ is called the \emph{ancestral configuration} of the gene tree at species-tree node $k$ under realization $R$. For example, in the realization of Fig.~\ref{configa}C, the ancestral configurations at the species-tree internal nodes are $C(f) = \{ a, b \}$, $C(g) = \{ d, e \}$, $C(h) = \{ c, f \}$, and $C(i) = \{ g, h \}$, where each gene lineage is identified by its immediate descendant node. In the realization of Fig.~\ref{configa}D, the ancestral configuration at the internal node of $S$ denoted by $h$ is $C(h) = \{ a, b, c \}$; at the other nodes, the ancestral configurations follow the previous case.

Let ${\bf R}(t)$ be the set of possible realizations of the gene-tree labeled topology $G = t$ in the matching species tree $S$. For a given node $k$ of $S$, by considering all possible realizations $R \in {\bf R}(t)$, we define the set
\begin{equation}\label{ciccio}
C_k = C_k(t) = \{ C(k,R) : R \in {\bf R}(t) \},
\end{equation}
with cardinality
\begin{equation}\label{cicci}
c_{k} = c_k(t) = |C_k|.
\end{equation}
Thus, $c_{k}$ counts ways that the gene lineages of $G$ can reach the time point right before node $k$ in $S$, when all realizations of the gene-tree labeled topology $G=t$ in $S$ are considered. For instance, taking $t$ as in Fig.~\ref{configa}A, $C_f = \{ \{ a, b \} \}$, $C_g = \{ \{ d, e \} \}$, 
\begin{small}$C_h = \{ \{a, b, c\}, \{ c, f \} \}$, and $C_i = \{ \{a, b, c, d, e\}, \{c, d, e, f\}, \{d, e, h\}, \{a, b, c, g\}, \{c, f, g\}, \{ g, h \} \}$\end{small}, for a total of 10 ancestral configurations. Note that from the definition of ancestral configuration, $\{k \} \notin C_k$. Indeed, gene lineages can coalesce to produce node $k$ of $G$ only in the part of the species tree above node $k$. For consistency with this observation, we set $c_k = 0$ if node $k$ is a leaf.

The set $C_k \cup \{ \{k \} \}$ can be viewed as the set of maximal antichains of the subtree $t^k$ of $t$ rooted at node $k$. In particular, if $r$ is the root of $t$, then $C_r \cup \{ \{r \} \}$ corresponds to the set of maximal antichains of $t$. An antichain of subtree $t^k$ is indeed a subset of its nodes---possibly including the leaves---whose elements are pairwise incomparable with respect to the descendant--ancestor order relation $\preceq$ defined in $t$. An ancestral configuration of $C_k$ is a ``maximal'' antichain of $t^k$ in the sense that it is not properly contained in any other antichain of $t^k$.

By summing the number $c_k$ for $k$ ranging over the set $N(t)$ of nodes of a labeled topology $t$, we find the total number of ancestral configurations of $t$, which we denote by
\begin{equation}\label{ciccitot}
c = c(t) = \sum_{k \in N(t)} c_{k}.
\end{equation}
Equivalently, $c+2|t|-1$ is the total number of maximal antichains across subtrees of $t$, including the $|t|$ leaves in $N(t)$ and for counts at each of the $|t|-1$ internal nodes, including as a maximal antichain the node itself.

For a gene-tree labeled topology $G=t$ and matching species tree $S$,
the total number of ancestral configurations $c(t)$ of $t$ is computed recursively by decomposing $t$ in its left and right root subtrees $t_L$ and $t_R$ (once we fix an embedding of $t$ into the plane). If $c_r(t)$ denotes the number of ancestral configurations at the root of $t$, then
\begin{eqnarray}\label{eqCtot}
c(t) & = & c(t_L) + c(t_R) + c_r(t) \\
\label{eqC}
c_r(t) & = & [c_r(t_L)+1][c_r(t_R)+1],
\end{eqnarray}
where $c(t) = c_r(t) = 0$ for $|t|=1$~\cite{DisantoAndRosenberg17}.

For example, suppose $t$ is the labeled topology of Fig.~\ref{configa}A, with internal nodes denoted as in Fig.~\ref{configa}B. Recalling that $t^k$ refers to the subtree of $t$ rooted at node $k$, by applying eqs.~(\ref{eqCtot}) and (\ref{eqC}), we find
\begin{eqnarray}\nonumber
c(t) &=& c(t^h) + c(t^g) + c_r(t) \\ \nonumber
& = & [c(t^f) + \bcancel{c(t^c)} + c_r(t^h)] +[\bcancel{c(t^d)} + \bcancel{c(t^e)} + c_r(t^g)] + [c_r(t^h)+1][c_r(t^g)+1] \\\nonumber
&=& c(t^f) + c_r(t^h) + c_r(t^g) + ([c_r(t^f)+1][\bcancel{c_r(t^c)}+1]  +1)( [\bcancel{c_r(t^d)}+1][\bcancel{c_r(t^e)}+1] +1) \\\nonumber
&=& [\bcancel{c(t^a)}+\bcancel{c(t^b)}+c_r(t^f)]+[c_r(t^f)+1][\bcancel{c_r(t^c)}+1] + [\bcancel{c_r(t^d)}+1][\bcancel{c_r(t^e)}+1] + [c_r(t^f)+2] 2 \\\nonumber
&=& 4 c_r(t^f) + 6 \\ \nonumber
&=& 4 [\bcancel{c_r(t^a)}+1][\bcancel{c_r(t^b)}+1] +6 \\ \nonumber
&=& 10. \nonumber
\end{eqnarray}

When the labeled topology $t$ has size $n$, the total number $c$ of ancestral configurations can be bounded by means of the number $c_r$ of root ancestral configurations as
\begin{equation}\label{bau}
c_r \leq c \leq (2n-1) c_r.
\end{equation}
Indeed, there are $|N(t)|=2n-1$ nodes in $t$ and, for every node $k$ of $t$, we have $c_k \leq c_r$.

Because $c$ and $c_r$ differ by a factor that is at most polynomial in the tree size $n$, they have the same exponential order when measured across tree families of increasing size. Based on this observation, the studies of \cite{DisantoEtAl20} and \cite{DisantoAndRosenberg17} of the asymptotic growth of the number of root ancestral configurations in random trees obtained the exponential order of the mean and variance of the total number of ancestral configurations in labeled topologies of size $n$ selected at random under the uniform and Yule distributions. In Section \ref{totconf}, we refine these results, obtaining full asymptotic distributions of the total number of ancestral configurations. We also study the correlation between the total number of ancestral configurations and the number of root ancestral configurations in random labeled topologies of increasing size.

\begin{figure}
\begin{center}
\includegraphics*[scale=0.40,trim=0 0 0 0]{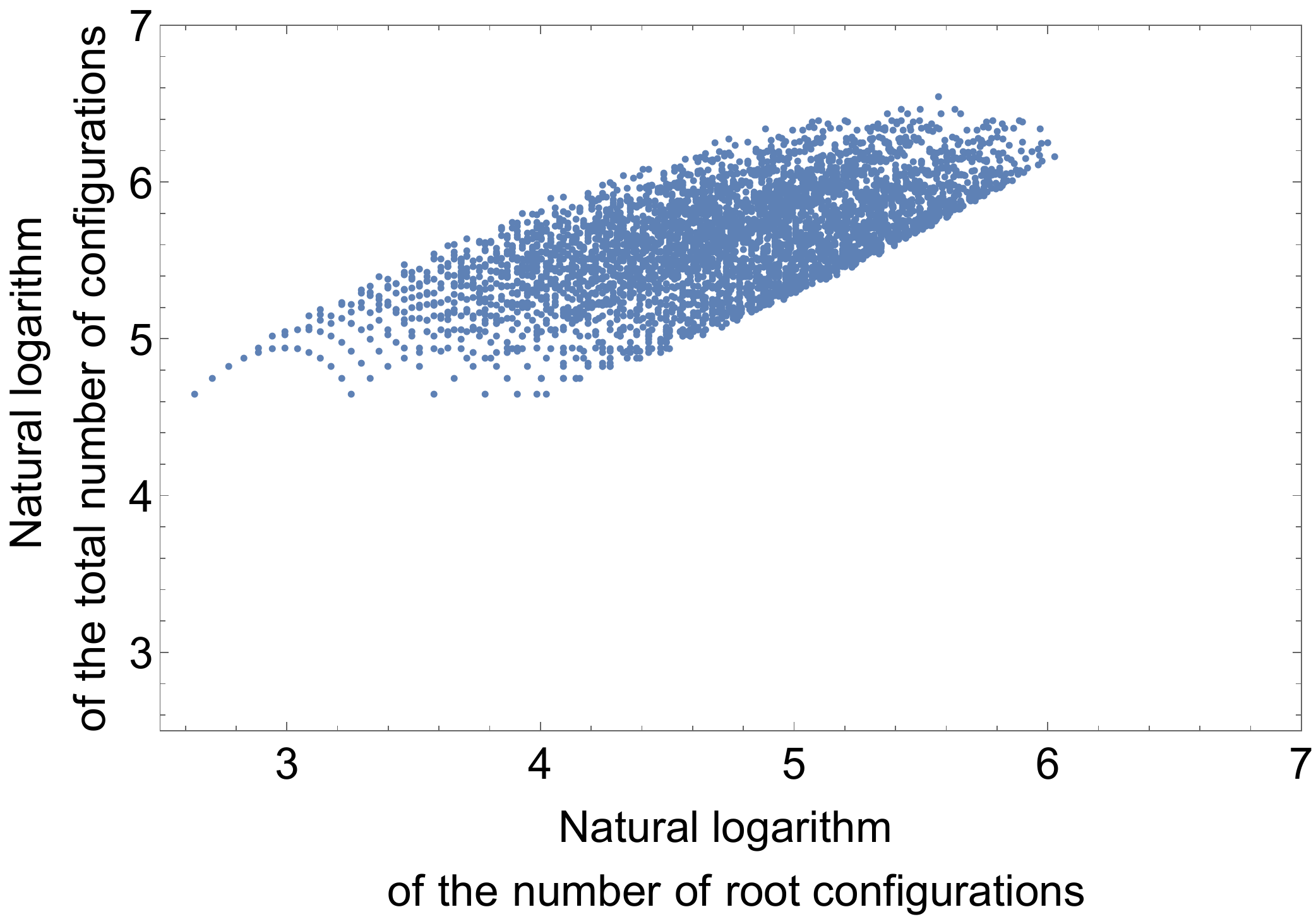}
\end{center}
\vspace{-0.4cm}
\caption{{\small Natural logarithms of the total number of ancestral configurations and the number of root ancestral configurations for representative labelings of each of the 4850 unlabeled topologies of size $n=15$ leaves. 
}} \label{totandroot}
\end{figure}

Fig.~\ref{totandroot} shows on a log scale the total number of ancestral configurations and the number of root ancestral configurations for representative labelings of each unlabeled topology of size $n=15$. The figure illustrates that the total number of ancestral configurations exceeds the number of root ancestral configurations. It also shows that the two quantities are positively correlated across trees.

\subsection{Ordered tree families and equivalent probability models of ancestral configurations}\label{models}

The definition in eq.~(\ref{ciccio}) of the set of ancestral configurations at a node of a labeled topology $t$ depends only on the shape of $t$. Ancestral configurations as well as the quantities in eqs.~(\ref{cicci}) and (\ref{ciccitot}) can be defined in the same way for many types of bifurcating rooted trees $t$ (e.g.~labeled, unlabeled, ordered, unordered). This section explains that probabilistic properties of the number of ancestral configurations considered over random labeled topologies can be equivalently analyzed over different tree families. In Sections \ref{fam1} and \ref{fam2}, we introduce the families of ordered unlabeled topologies and ordered unlabeled histories. Next, in Section \ref{eqfam}, we recall some equivalence results of \cite{DisantoEtAl20} on the distribution of the number of ancestral configurations. In particular, Lemma \ref{lemcorr} states that the number of ancestral configurations has the same distribution when considered over uniformly distributed labeled topologies and over uniformly distributed ordered unlabeled topologies of the same size. The equivalence extends to the distribution of the number of ancestral configurations over Yule-distributed labeled topologies and uniformly distributed ordered unlabeled histories. We conclude the section with Lemma \ref{rectot}, a preliminary result to the calculations of Section \ref{totconf}.

\subsubsection{Ordered unlabeled topologies}\label{fam1}

An \emph{ordered} unlabeled topology is a binary rooted plane tree, that is, an unlabeled topology $t$ equipped with a left--right orientation of the subtrees descending from its internal nodes. The tree \emph{shape} of an ordered unlabeled topology is the underlying unordered unlabeled topology. Each ordered unlabeled topology is an embedding of its shape into the plane.

\begin{figure}
\begin{center}
\includegraphics*[scale=0.56,trim=0 0 0 0]{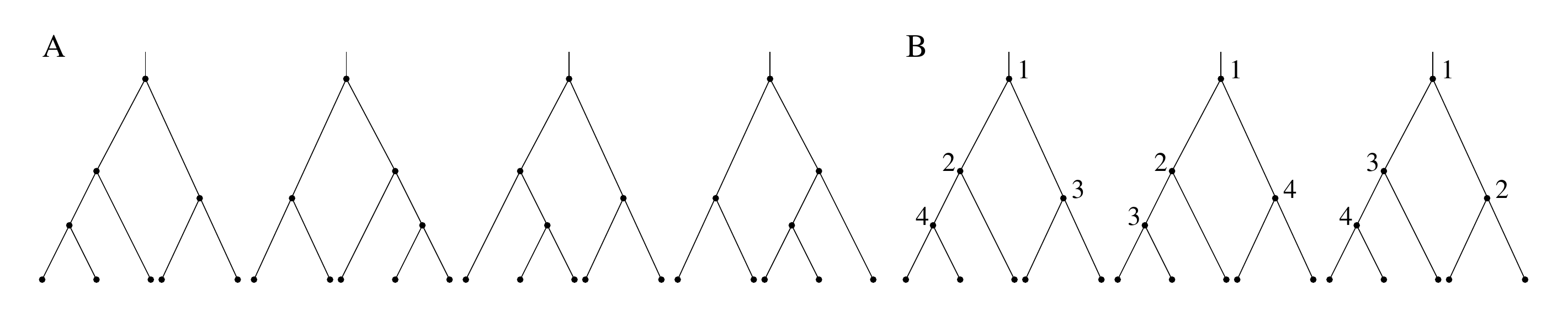}
\end{center}
\vspace{-.7cm}
\caption{{\small
Ordered unlabeled topologies and ordered unlabeled histories. {\bf(A)} The four possible ordered unlabeled topologies whose underlying unordered unlabeled topology is  $(((\bullet,\bullet),\bullet),(\bullet,\bullet))$. {\bf(B)} The three possible ordered unlabeled histories whose underlying ordered unlabeled topology matches the leftmost shape in (A).}} \label{orient}
\end{figure}

Fig.~\ref{orient}A depicts the four ordered unlabeled topologies with shape $(((\bullet,\bullet),\bullet),(\bullet,\bullet))$. Denoting by ${\rm out}(t)$ the number of \underline{o}rdered \underline{u}nlabeled \underline{t}opologies with shape $t$, eq.~(23) of \cite{DisantoEtAl20} gives
\begin{equation}\label{outeq}
{\rm out}(t) = 2 \, {\rm out}(t_L) \, {\rm out}(t_R) \, \frac{1}{1 + \delta_{t_L = t_R}},
\end{equation}
where ${\rm out}(t) = 1$ if $|t| = 1$.

Ordered unlabeled topologies are also called \emph{Catalan trees} as they are enumerated, with respect to the size $n$, by
the $(n-1)$-th Catalan number $C_{n-1}$ (\cite{Stanley99}, Exercise 6.19d), where
\begin{equation}\label{cata}
C_n = \frac{1}{n+1} {{2n}\choose{n}} \sim \frac{4^n}{\sqrt{\pi n^3}}.
\end{equation}
The generating function of the sequence $C_n$ is
\begin{equation}\label{catag}
C(z) = \sum_{n=0}^\infty C_n z^n = \frac{1-\sqrt{1-4z}}{2z}.
\end{equation}
$C(z)$ has singular expansion $C(z) \stackrel{z \rightarrow \alpha}{\sim} 2 - 2\sqrt{1-4z}$ at its dominant singularity $\alpha=\frac{1}{4}$, as can be seen by 
setting $z=\frac{1}{4}$ 
in the denominator; 
 by eq.~(\ref{pillo2}), we obtain the asymptotic expression in eq.~(\ref{cata}).

A decomposition provides a useful formula for the probability that the left and right root subtrees of a uniformly selected ordered unlabeled topology of $n$ leaves have sizes $j$ and $n-j$, respectively ($1 \leq j \leq n$). Each ordered unlabeled topology $t$ of $n \geq 2$ leaves results from the following recursive construction: (a) take two ordered unlabeled topologies $t_{L}$ and $t_R$ of sizes $j$ and $n-j$, respectively, and (b) append $t_{L}$ and $t_R$ to the left and right, respectively, of a common root node. For example, the leftmost ordered unlabeled topology of Fig.~\ref{orient}A is obtained by appending $t_{L} = ((\bullet,\bullet),\bullet)$ to the left and $t_R = (\bullet,\bullet)$ to the right of the shared root. For the third ordered unlabeled topology of Fig.~\ref{orient}A, we take instead $t_{L} = (\bullet,(\bullet,\bullet)))$ and $t_R = (\bullet,\bullet)$. Because $C_{j-1}$ possible choices exist for $t_{L}$ and $C_{n-1-j}$ exist for $t_{R}$, the probability that a uniformly distributed ordered unlabeled topology $t$ of size $n$ has a left root subtree $t_{L}$ with size $j$ and a right root subtree $t_R$ with size $n-j$ is
\begin{equation}\label{trilli}
\mathbb{P}[|t_{L}| = j \, \& \, |t_{R}| = n-j] = \frac{C_{j-1} \, C_{n-1-j}}{C_{n-1}}.
\end{equation}

\subsubsection{Ordered unlabeled histories}\label{fam2}

An ordered unlabeled \emph{history} of size $n$ is a plane embedding of an unlabeled topology of $n$ leaves whose internal nodes are bijectively labeled by the integers from the interval $[1,n-1]$ in such a way that each non-root internal node has a larger label than its parent node (Fig.~\ref{orient}B). From a biological standpoint, the labels at the internal nodes define a temporal ordering of the coalescent events in the history.

In the language of computer science, ordered unlabeled histories---with leaves and their incident edges stripped away---correspond to the so-called \emph{increasing binary trees} (Example II.17 of \cite{FlajoletAndSedgewick09}), with the term ``increasing'' referring to the labels of the nodes that increase along any path from the root to the leaves of the tree. To specify the linear ordering of the internal nodes, we write the Newick format of an ordered unlabeled history by adding as a subscript next to a closed parenthesis the label of the corresponding internal node. For instance, $(((\bullet,\bullet)_4,\bullet)_2,(\bullet,\bullet)_3)_1$ indicates the first ordered unlabeled history depicted in Fig.~\ref{orient}B.

The \emph{shape} of an ordered unlabeled history is the underlying unordered unlabeled topology obtained by removing labels at internal nodes and ignoring left--right orientation. With the same notation used in eqs.~(\ref{lteq}) and (\ref{outeq}), the number ${\rm ouh}(t)$ of \underline{o}rdered \underline{u}nlabeled \underline{h}istories with tree shape $t$ is calculated recursively as
\begin{equation}\label{ouheq}
{\rm ouh}(t) = 2 \, {\rm ouh}(t_L) \, {\rm ouh}(t_R) \, {{|t|-2}\choose{|t_L|-1}} \, \frac{1}{1 + \delta_{t_L = t_R}},
\end{equation}
where ${\rm ouh}(t) = 1$ if $|t|=1$. Each ordered unlabeled history with shape $t$ is constructed by appending two ordered unlabeled histories $h_1$ and $h_2$ with shapes $t_L$ and $t_R$, respectively, to the left and right of a common root node, while choosing one of ${{|t|-2}\choose{|t_L|-1}}$ possibilities for merging the linear ordering of the internal nodes of $h_1$ with that of the internal nodes of $h_2$. The factor $1/(1 + \delta_{t_L = t_R})$ accounts for possible symmetries in this process.

The set of all possible ordered unlabeled histories of size $n$ is enumerated by $F_{n-1}$ (\cite{Steel16}, p.~47), where
\begin{equation}\label{fact}
F_n = n!.
\end{equation}
A formula analogous to eq.~(\ref{trilli}) can be found for ordered unlabeled histories by extending the recursive construction of ordered unlabeled topologies. To construct an ordered unlabeled history $t$ of size $n$, we do the following: (a) take two ordered unlabeled histories $t_{L}$ and $t_R$ of sizes $j$ and $n-j$, respectively, (b) append $t_{L}$ and $t_R$ to the left and to the right, respectively, of a shared root node, and (c) merge the linear ordering of the internal nodes of $t_{L}$ with the linear ordering of the internal nodes of $t_{R}$ to define a linear ordering of the internal nodes of $t$.

For example, the leftmost ordered unlabeled history of Fig.~\ref{orient}B is obtained by appending $t_{L} = ((\bullet,\bullet)_2,\bullet)_1$ to the left and $t_R = (\bullet,\bullet)_1$ to the right of the same root, and then merging the orderings of the internal nodes of $t_{L}$ and $t_R$ by putting the root node of $t_R$ between the internal nodes of $t_{L}$. Because there are $F_{j-1}$ choices for $t_{L}$, $F_{n-1-j}$ choices for $t_{R}$, and ${{n-2}\choose{j-1}}$ ways to merge the ordering of the $j-1$ internal nodes of $t_{L}$ with the ordering of the $n-1-j$ internal nodes of $t_R$, the probability that a uniformly distributed ordered unlabeled history $t$ of size $n$ has its left root subtree $t_{L}$ of size $j$ and its right root subtree $t_R$ of size $n-j$ is
\begin{equation}\label{trilli2}
\mathbb{P}[|t_{L}| = j \, \& \, |t_{R}| = n - j] = \frac{F_{j-1} \, F_{n-1-j} \, {{n-2}\choose{j-1}}}{F_{n-1}} = \frac{1}{n-1}.
\end{equation}

\subsubsection{Equivalent models for ancestral configurations} \label{eqfam}

\begin{table}\caption{Induced probabilities of unlabeled topologies of size 5.}
\vspace{-.5cm}
\begin{center}\label{induc}
\fontsize{9}{11}\selectfont
\begin{tabular}{||c||ccc|ccc||}
\hline \hline
& Probability induced & & Probability induced & Probability induced & & Probability induced\\
Unlabeled & by uniform ordered   & & by uniform        & by uniform ordered  & & by Yule\\
topology  & unlabeled topologies & & labeled toplogies & unlabeled histories & & labeled topologies \\
& & & & & & \\
$t$ & $\frac{{\rm out}(t) }{ C_4 } $ && $\frac{{\rm lab}(t) }{ |L_5| }$ & $\frac{{\rm ouh}(t) }{ F_4 }$ && ${\rm lab}(t) \times \mathbb{P}_Y[t]$ \\[2ex] \hline
& & & & & & \\
\includegraphics[scale=.12,bb=200 200 210 250]{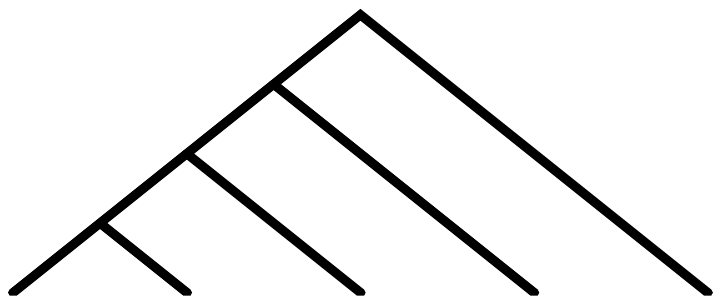} & $\frac{8}{14} = \frac{4}{7}$ && $\frac{60}{105} = \frac{4}{7}$ & $\frac{8}{24} = \frac{1}{3} $ && $60 \times \frac{1}{180} = \frac{1}{3} $\\[2ex]
\includegraphics[scale=.12,bb=200 200 210 250]{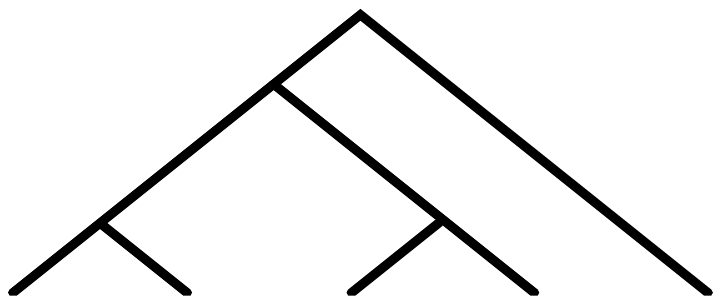} & $\frac{2}{14} = \frac{1}{7}$ && $\frac{15}{105} = \frac{1}{7}$ & $\frac{4}{24} = \frac{1}{6} $ && $15 \times \frac{1}{90} = \frac{1}{6} $\\[2ex]
\includegraphics[scale=.12,bb=200 200 210 250]{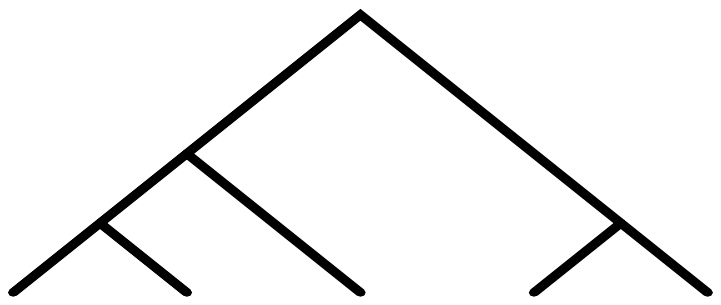} & $\frac{4}{14} = \frac{2}{7}$ && $\frac{30}{105} = \frac{2}{7} $ & $\frac{12}{24} = \frac{1}{2} $ && $30 \times \frac{1}{60} = \frac{1}{2} $\\[2ex]
\hline \hline
\end{tabular}
\end{center}
\vspace{-.1cm}
{\small{For each unlabeled topology $t$ of size $5$, the probability of $t$ induced by the uniform distribution over ordered unlabeled topologies is the ratio of ${\rm out}(t)$ to the total number $C_4=14$ of ordered unlabeled topologies. Similar calculations appear for uniform labeled topologies ($|L_5|=105$), ordered unlabeled histories ($F_4=24$), and Yule labeled topologies. Quantities ${\rm out}(t)$, ${\rm lab}(t)$, and ${\rm ouh}(t)$ are recursively computed from eqs.~(\ref{outeq}), (\ref{lteq}), and (\ref{ouheq}), respectively. The probability of $t$ induced by the Yule distribution over labeled topologies is obtained by multiplying ${\rm lab}(t)$ by the Yule probability in eq.~(\ref{pyule}) of each labeled topology with $t$ as underlying unlabeled topology. The second and third columns agree, as do the fourth and fifth.}}
\end{table}

We previously noticed~\cite{DisantoEtAl20} that ordered unlabeled topologies and ordered unlabeled histories can be used to study the number of ancestral configurations of uniformly and Yule-distributed labeled topologies, respectively. We observed that the number of ancestral configurations of a given tree structure depends only on the underlying unlabeled topology. Second, as shown in the proofs of Lemmas 1, 2 and 3 of \cite{DisantoEtAl20}, the uniform distribution over the set of ordered unlabeled topologies of size $n$ and the uniform distribution over the set of labeled topologies of size $n$ induce the same distribution over the set of underlying unlabeled topologies of size $n$; the uniform distribution over ordered unlabeled histories of size $n$ and the Yule distribution over labeled topologies of size $n$ induce the same distribution over the unlabeled topologies of size $n$. In other words, for each unlabeled topology $t$, the sum of the probabilities of the uniformly distributed ordered unlabeled topologies (resp.~histories) having the shape of $t$ equals the sum of the probabilities of the uniformly (resp.~Yule) distributed labeled topologies with tree shape $t$. These two facts yield the next lemma. Table \ref{induc} shows the case $n=5$.
\begin{lemm}\label{lemcorr}
The distribution of the number of ancestral configurations over uniformly (resp.~Yule) distributed labeled topologies of size $n$ is the distribution of the number of ancestral configurations over uniformly distributed ordered unlabeled topologies (resp.~histories) of size $n$.
\end{lemm}
From this lemma, probabilistic properties of the number of ancestral configurations of uniformly and Yule-distributed labeled topologies can be equivalently studied over uniformly distributed ordered unlabeled topologies and uniformly distributed ordered unlabeled histories, respectively.

To use these equivalences, we require the following lemma.
\begin{lemm}
\label{rectot}
Let $R_n$ and $T_n$ be the random variables that represent the number of root ancestral configurations and the total number of ancestral configurations in a random ordered unlabeled topology (resp.~history) of size $n$ selected under the uniform distribution. Equivalently, by Lemma \ref{lemcorr}, $R_n$ and $T_n$ represent the numbers of root ancestral configurations and the total number of ancestral configurations in a random labeled topology of size $n$ selected under the uniform (resp.~Yule) distribution. Then we have $R_1=T_1=0$, and for $n\geq 2$,
\begin{eqnarray}
\label{lemm2line1}
T_n&\stackrel{d}{=} & T_{I_n}+T_{n-I_n}^{*}+R_n,\\
\label{lemm2line2}
R_n&\stackrel{d}{=} & R_{I_n} \, R^{*}_{n-I_n} + R_{I_n} + R^{*}_{n-I_n} + 1,
\end{eqnarray}
where $I_n$ is distributed over the interval $[1,n-1]$ with probability $\mathbb{P}[I_n=j]={C_{j-1} \, C_{n-1-j}} / {C_{n-1}}$ (resp. $\mathbb{P}[I_n=j]=\frac{1}{n-1}$), $R_j^{*}$ and $T_j^{*}$ are independent copies of $R_j$ and $T_j$, respectively, for each $j \in [1,n-1]$, and both $R_j$ and $R_j^{*}$ as well as $T_j$ and $T_j^{*}$ are independent of $I_j$ for $j \in [1,n-1]$.
\end{lemm}
\noindent \emph{Proof.} The distributional recurrences follow directly from eqs.~(\ref{eqCtot}) and (\ref{eqC}). $\mathbb{P}[I_n=j]$ follows eqs.~(\ref{trilli}) and (\ref{trilli2}), giving the probability that the left root subtree of an ordered unlabeled topology or history of $n$ taxa selected uniformly at random has size $I_n = j$. $\Box$

\subsection{Known results on the distribution of ancestral configurations}
\label{distrib}

For the random variables $R_n$ and $T_n$, the asymptotic behavior of the moments $\mathbb{E}[R_n]$, $\mathbb{E}[R_n^2]$, and $\mathbb{E}[T_n]$ and variances $\mathbb{V}[R_n]$ and $ \mathbb{V}[T_n]$ were studied under the uniform model of labeled topologies by \cite{DisantoAndRosenberg17} (Propositions 5 and 6), and under the Yule model by \cite{DisantoEtAl20} (Propositions 5.4 and 5.5):
\begin{equation}\label{meandis-propyu}
    \mathbb{E}[T_n] \bowtie \mathbb{E}[R_n] \sim
    \begin{cases}
       \sqrt{\frac{3}{2}} \left(\frac{4}{3}\right)^n, & \text{Uniform model,} \\
      \left(\frac{1}{ 1-e^{-2\pi\sqrt{3}/9} } \right)^n, & \text{Yule model,}
    \end{cases}
  \end{equation}
\begin{equation}\label{vardis-propyuv}
    \mathbb{V}[T_n] \bowtie \mathbb{V}[R_n] \sim \mathbb{E}[R_n^2] \sim
    \begin{cases}
      \sqrt{\frac{7(11 - \sqrt{2})}{34}}\left[\frac{4}{7 (8\sqrt{2} - 11) }\right]^n, & \text{Uniform model,} \\
      (2.0449954971\ldots)^n, & \text{Yule model.}
    \end{cases}
  \end{equation}
The exponential order $2.0449954971\ldots$ of $\mathbb{V}[T_n]$ under the Yule model was approximated by a numerical procedure described in the Appendix of \cite{DisantoEtAl20}.

Under both the uniform and Yule models, the logarithm of the number of root configurations of a randomly selected labeled topology of size $n$ was shown to asymptotically follow a normal distribution.
Propositions 4.1 and 5.2 of \cite{DisantoEtAl20} state that the rescaled random variable $$\frac{\log R_n - \mathbb{E}[\log R_n]}{\sqrt{\mathbb{V}[\log R_n]}}$$ converges to a standard normal distribution, where
$\mathbb{E}[\log R_n ] \sim \mu n$, $\mathbb{V}[\log R_n ] \sim \sigma^2 n$, and
\begin{equation}\label{lognormalconst}
(\mu, \sigma^2) \approx
\begin{cases}
  (0.272, 0.034), & \text{Uniform model,} \\
  (0.351, 0.008), & \text{Yule model}.
\end{cases}
\end{equation}

\section{Distributional properties of the total number of ancestral configurations}\label{totconf}

Previous results on ancestral configurations focused on the number of root configurations of labeled topologies selected under the uniform and Yule distributions. We now study the random \emph{total number of ancestral configurations} under the same two probability models. In particular, we determine the asymptotic growth of its mean and variance. In agreement with eqs.~(\ref{meandis-propyu}) and (\ref{vardis-propyuv}), we find that the mean and variance of the total number of configurations differ from the mean and variance of the number of root configurations only in their subexponential terms, which turn out to be constants. Moreover, we find that, as is true of the number of root configurations, the total number of configurations follows an asymptotically lognormal distribution.

\subsection{Uniform ordered unlabeled topologies and uniform labeled topologies}

By Lemma \ref{lemcorr}, the distribution of the number of ancestral configurations over random labeled topologies of size $n$ selected uniformly at random is the distribution of the number of ancestral configurations over uniformly distributed ordered unlabeled topologies of size $n$. We use this equivalence to derive the results of this section, denoting by $R_n$ and $T_n$, respectively, the number of root ancestral configurations and the total number of ancestral configurations in a random ordered unlabeled topology of size $n$ selected under the uniform distribution.

Our first proposition uses the system of distributional recurrences of Lemma \ref{rectot} to determine the asymptotic behavior of the mean of $T_n$.
\begin{prop}\label{meantotuni}
The mean total number of ancestral configurations in an ordered unlabeled topology of size $n$ selected uniformly at random satisfies the asymptotic relation $\mathbb{E}[T_n] \sim 2 \mathbb{E}[R_n] \sim \sqrt{6} (4/3)^n$.
\end{prop}
\noindent \emph{Proof.}
By eq.~(\ref{lemm2line2}) in Lemma \ref{rectot} coupled with $\mathbb{E}[ R_{I_n} \, R^{*}_{n-I_n}] = \sum_{j=1}^{n-1} \mathbb{P}[I_n = j] \, \mathbb{E}[ R_{j} \, R^{*}_{n-j}] = \sum_{j=1}^{n-1} \mathbb{P}[I_n = j] \, \mathbb{E}[ R_{j}] \, \mathbb{E}[R^{*}_{n-j}]$, we find that for $n \geq 1$, the expectation of $R_n$ satisfies
\begin{equation}\label{ern1}
C_{n-1} {\mathbb E}[R_n]=\sum_{j=1}^{n-1}C_{j-1} \, C_{n-1-j} \, \Big({\mathbb E}[R_j] \, {\mathbb E}[R_{n-j}]+{\mathbb E}[R_j]+{\mathbb E}[R_{n-j}]+1\Big),
\end{equation}
which holds also for $n=1$ as $\mathbb{E}[ R_{1}] = 0$. Similarly, with $\mathbb{E}[ T_{1}] = 0$, for $n \geq 1$, eq.~(\ref{lemm2line1}) in Lemma \ref{rectot} gives
\begin{equation}\label{etn1}
C_{n-1} {\mathbb E}[T_n]=2 \left(\sum_{j=1}^{n-1}C_{j-1}C_{n-1-j}{\mathbb E}[T_j] \right)+ C_{n-1}{\mathbb E}[R_n].
\end{equation}

Define the generating functions
\begin{eqnarray} \label{rz}
R(z) & \equiv & \sum_{n=1}^\infty C_{n-1}{\mathbb E}[R_n]z^n \\
\label{tz}
T(z) & \equiv & \sum_{n=1}^\infty C_{n-1}{\mathbb E}[T_n]z^n,
\end{eqnarray}
whose coefficients $C_{n-1}{\mathbb E}[R_n] = C_{n-1} \sum_{i=0}^\infty {|\{ t : c_r(t) = i \}|}i/{C_{n-1}}$ and $C_{n-1}{\mathbb E}[T_n] = C_{n-1} \sum_{i=0}^\infty {|\{ t : c(t) = i \}|}i/{C_{n-1}}$ give respectively the sum of the number of root configurations and the sum of the total number of configurations over all ordered unlabeled topologies $t$ of $n$ taxa.

The recurrences in eqs.~(\ref{ern1}) and (\ref{etn1}) translate into a system of equations for $R(z)$ and $T(z)$:
\[
S_1 \equiv \left\{
\begin{array}{ll}
R(z) = R(z)^2 + 2zC(z) \, R(z) + z^2C(z)^2 \\
T(z) = 2zC(z) \, T(z) + R(z),
\end{array}
\right.
\]
where $C(z)$ is the Catalan generating function (eq.~\ref{catag}). Indeed, multiplying eq.~(\ref{ern1}) by $z^n$, we have
\begin{eqnarray}\nonumber
C_{n-1} \mathbb{E}[R_n] z^n &=& \sum_{j=1}^{n-1} C_{j-1} {\mathbb E}[R_j] z^j \cdot C_{n-1-j} {\mathbb E}[R_{n-j}] z^{n-j}
+ z \sum_{j=1}^{n-1} C_{j-1} {\mathbb E}[R_j] z^j \cdot C_{n-1-j} z^{n-1-j} \\\nonumber
&& + z \sum_{j=1}^{n-1} C_{j-1} z^{j-1} \cdot C_{n-1-j} {\mathbb E}[R_{n-j}] z^{n-j}
+ z^2 \sum_{j=1}^{n-1} C_{j-1} z^{j-1} \cdot C_{n-1-j} z^{n-1-j}.
\end{eqnarray}
The first equation of $S_1$ is obtained by summing over $n \geq 1$ and simplifying:
\begin{eqnarray}\nonumber
R(z) &=& \sum_{n=1}^\infty C_{n-1} \mathbb{E}[R_n] z^n = \sum_{n=1}^\infty \sum_{j=1}^{n-1} C_{j-1} {\mathbb E}[R_j] z^j \cdot C_{n-1-j} {\mathbb E}[R_{n-j}] z^{n-j}
+ z \sum_{n=1}^\infty \sum_{j=1}^{n-1} C_{j-1} {\mathbb E}[R_j] z^j \cdot C_{n-1-j} z^{n-1-j}  \\\nonumber
&& + z \sum_{n=1}^\infty \sum_{j=1}^{n-1} C_{j-1} z^{j-1} \cdot C_{n-1-j} {\mathbb E}[R_{n-j}] z^{n-j}
+ z^2 \sum_{n=1}^\infty \sum_{j=1}^{n-1} C_{j-1} z^{j-1} \cdot C_{n-1-j} z^{n-1-j}.
\end{eqnarray}
Similarly, from eq.~(\ref{etn1}), we obtain the second equation of $S_1$:
\begin{eqnarray}\nonumber
T(z) = \sum_{n=1}^\infty C_{n-1} {\mathbb E}[T_n] z^n &= &2 z \sum_{n=1}^\infty \sum_{j=1}^{n-1}C_{j-1}\mathbb{E}[T_j] z^j \cdot C_{n-1-j} z^{n-1-j} + \sum_{n=1}^\infty C_{n-1} {\mathbb E}[R_n] z^n.
\end{eqnarray}

Solving system $S_1$ for $T(z)$ yields
\begin{equation}\label{eqtt}
T(z)=\frac{R(z)}{1-2zC(z)}=\frac{\sqrt{1-4z}-\sqrt{2\sqrt{1-4z}-1}}{2\sqrt{1-4z}},
\end{equation}
which has dominant singularity $\alpha \equiv \frac{3}{16}$---the root of $2\sqrt{1-4z}-1$. The singular expansion is
$$T(z) \stackrel{z \rightarrow \alpha}{\sim} k_1 - \sqrt{\frac{3}{2}} \cdot \sqrt{1-\frac{16 z}{3}},$$
for a certain constant $k_1$. Eq.~(\ref{pillo2}) thus yields
$$[z^n] T(z) \sim \sqrt{\frac{3}{2}} \frac{(16/3)^n}{2 \sqrt{\pi n^3}}.$$
By using the fact that $C_{n-1} \sim {4^{n-1}}/{\sqrt{\pi n^3}}$ (eq.~\ref{cata}), we obtain
\[
{\mathbb E}[T_n] = \frac{[z^n]T(z)}{C_{n-1}} \sim \frac{\sqrt{\frac{3}{2}} \frac{(16/3)^n}{2 \sqrt{\pi n^3}}}{\frac{4^{n-1}}{\sqrt{\pi n^3}}} = \sqrt{6}\left(\frac{4}{3}\right)^n,
\] which is twice the asymptotic value of ${\mathbb E}[R_n]$ given for the uniform case in eq.~(\ref{meandis-propyu}). $\Box$

\begin{figure}
\begin{center}
\includegraphics*[scale=0.40,trim=0 0 0 0]{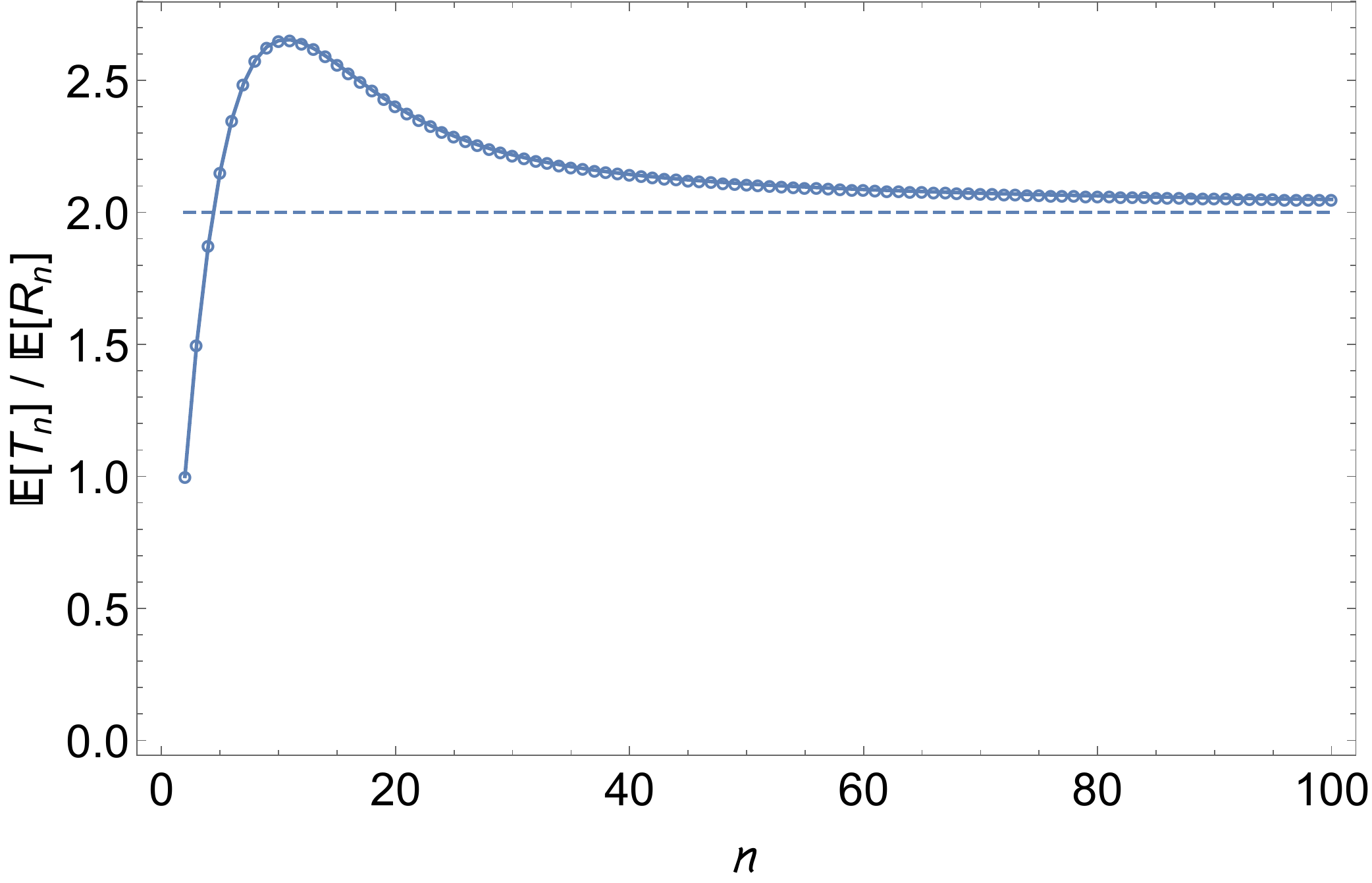}
\end{center}
\vspace{-.7cm}
\caption{{\small Ratio of the mean total number $\mathbb{E}[T_n]$ of ancestral configurations and mean number $\mathbb{E}[R_n]$ of root configurations for uniformly distributed ordered unlabeled topologies (or uniformly distributed labeled topologies) of size $2 \leq n \leq 100$. Values of $\mathbb{E}[R_n]$ and $\mathbb{E}[T_n]$ are computed from the recurrences in eqs.~(\ref{ern1}) and (\ref{etn1}), respectively.}}
\label{ratiounif}
\end{figure}

\vskip .2cm
Fig.~\ref{ratiounif} plots the exact ratio $\mathbb{E}[T_n] / \mathbb{E}[R_n]$ with increasing $n$. In agreement with Proposition \ref{meantotuni}, the ratio $\mathbb{E}[T_n] / \mathbb{E}[R_n]$ approaches 2 as $n$ increases.

We now consider the variance ${\mathbb V}[T_n]$ of the total number of ancestral configurations and its correlation coefficient $\rho [T_n, R_n]$ with the number of root configurations in uniformly distributed ordered unlabeled topologies of fixed size $n$. The next lemma provides a series of distributional recurrences.
\begin{lemm}\label{4rec}
Consider the random variables $\tilde{R}_n \equiv R_n+1$ and $T_n$. We have $\tilde{R}_1=1$, $T_1=0$, and for $n \geq 2$,
\begin{eqnarray}
\label{sec-r}
\tilde{R}_n &\stackrel{d}{=} & \tilde{R}_{I_n}\tilde{R}^{*}_{n-I_n}+1, \\
\label{sec-rn}
(\tilde{R}_n)^2 & \stackrel{d}{=} &(\tilde{R}_{I_n})^2(\tilde{R}^{*}_{n-I_n})^2+2\tilde{R}_{I_n}\tilde{R}^{*}_{n-I_n}+1, \\
\label{sec-rnn}
T_n\tilde{R}_n &\stackrel{d}{=} & T_{I_n}\tilde{R}_{I_n}\tilde{R}^*_{n-I_n}+T^{*}_{n-I_n}\tilde{R}^{*}_{n-I_n}\tilde{R}_{I_n}+T_{I_n}+T^{*}_{n-I_n}+(\tilde{R}_n)^2-\tilde{R}_n, \\
\label{sec-rnnn}
(T_n)^2&\stackrel{d}{=} & (T_{I_n})^2+(T^{*}_{n-I_n})^2+2T_{I_n}T^{*}_{n-I_n}+2T_nR_n-(R_n)^2,
\end{eqnarray}
where $I_n$ is distributed over the interval $[1,n-1]$ with probability $\mathbb{P}[I_n=j]={C_{j-1}C_{n-1-j}}/{C_{n-1}}$, $R_j^{*}$, $\tilde{R}^*_j$, and $T_j^{*}$ are independent copies of $R_j$, $\tilde{R}_j$ and $T_j$, respectively, for every $j \in [1,n-1]$, and both $R_j$ and $R_j^{*}$ as well as $\tilde{R}_j$, $\tilde{R}^*_j$, $T_j,$ and $T_j^{*}$ are independent of $I_j$ for $j \in [1,n-1]$.
\end{lemm}
\noindent \emph{Proof.} Eq.~(\ref{sec-r}) follows directly from  eq.~(\ref{lemm2line2}) in Lemma \ref{rectot}. Eq.~(\ref{sec-rn}) is obtained by squaring eq.~(\ref{sec-r}). For eq.~(\ref{sec-rnn}), eq.~(\ref{lemm2line1}) in Lemma \ref{rectot} with eq.~(\ref{sec-r}) yields
$$T_n\tilde{R}_n = (T_{I_n}+T_{n-I_n}^{*}+R_n)\tilde{R}_n = (T_{I_n}+T_{n-I_n}^{*})\tilde{R}_n + R_n\tilde{R}_n = (T_{I_n}+T_{n-I_n}^{*})( \tilde{R}_{I_n}\tilde{R}^{*}_{n-I_n}+1) + (\tilde{R}_n-1) \tilde{R}_n.$$
Finally, by squaring eq.~(\ref{lemm2line1}) in Lemma \ref{rectot}, we obtain
\begin{eqnarray*}
(T_n)^2&\stackrel{d}{=} &(T_{I_n})^2+(T^{*}_{n-I_n})^2+(R_n)^2+2T_{I_n}T^{*}_{n-I_n}+2T_{I_n}R_n+2T^{*}_{n-I_n}R_n\\
& \stackrel{d}{=}& (T_{I_n})^2+(T^{*}_{n-I_n})^2+ 2T_{I_n}T^{*}_{n-I_n} + 2(R_n)^2 + 2T_{I_n}R_n+2T^{*}_{n-I_n}R_n - (R_n)^2 \\
& \stackrel{d}{=}& (T_{I_n})^2+(T^{*}_{n-I_n})^2+ 2T_{I_n}T^{*}_{n-I_n} + 2R_n ( R_n + T_{I_n} + T^{*}_{n-I_n} ) - (R_n)^2,
\end{eqnarray*}
which gives eq.~(\ref{sec-rnnn}), because $R_n + T_{I_n} + T^{*}_{n-I_n} = T_n$ again by eq.~(\ref{lemm2line1}). $\Box$

\vskip .2cm
To proceed with the asymptotic analysis of the variance ${\mathbb V}[T_n]$ and correlation $\rho [T_n, R_n]$, we now determine the asymptotic behavior of expectations ${\mathbb E}[T_n^2]$ and ${\mathbb E}[T_n R_n]$. We define the following generating functions:
\begin{eqnarray}\label{erretilde}
\tilde{R}(z) &\equiv & \sum_{n=1}^\infty C_{n-1}{\mathbb E}[\tilde{R}_n]z^n = \sum_{n=1}^\infty C_{n-1}{\mathbb E}[R_n]z^n + \sum_{n=1}^\infty C_{n-1}z^n = R(z)+zC(z),\\
\tilde{S}(z) &\equiv &\sum_{n=1}^\infty C_{n-1}{\mathbb E}[\tilde{R}_n^2]z^n, \\ \nonumber
S(z) &\equiv &\sum_{n=1}^\infty C_{n-1}{\mathbb E}[R_n^2]z^n = \sum_{n=1}^\infty C_{n-1}{\mathbb E}[\tilde{R}_n^2 - 2 \tilde{R}_n +1]z^n = \sum_{n=1}^\infty C_{n-1}{\mathbb E}[\tilde{R}_n^2 - 2 (R_n + 1) +1]z^n \\ \label{esse}
&=& \sum_{n=1}^\infty C_{n-1}{\mathbb E}[\tilde{R}_n^2]z^n -2 \sum_{n=1}^\infty C_{n-1}{\mathbb E}[R_n]z^n - \sum_{n=1}^\infty C_{n-1}z^n = \tilde{S}(z)-2R(z)-zC(z), \\
\tilde{V}(z) &\equiv &\sum_{n=1}^\infty C_{n-1}{\mathbb E}[T_n\tilde{R}_n]z^n,\\ \nonumber
V(z) &\equiv & \sum_{n=1}^\infty C_{n-1}{\mathbb E}[T_nR_n]z^n = \sum_{n=1}^\infty C_{n-1}{\mathbb E}[T_n(\tilde{R}_n-1)]z^n = \sum_{n=1}^\infty C_{n-1}{\mathbb E}[T_n\tilde{R}_n]z^n - \sum_{n=1}^\infty C_{n-1}{\mathbb E}[T_n]z^n \\ \label{vu}
&=& \tilde{V}(z)-T(z), \\
U(z) & \equiv & \sum_{n=1}^\infty C_{n-1}{\mathbb E}[T_n^2]z^n.
\end{eqnarray}
Here, $R(z)$ is from eq.~(\ref{rz}), $C(z)$ is given in eq.~(\ref{catag}), and $T(z)$ is given in eq.~(\ref{eqtt}).

The distributional recurrences of Lemma \ref{4rec} determine recurrences for the expectations ${\mathbb E}[\tilde{R}_n^2]$, ${\mathbb E}[T_n\tilde{R}_n]$, and ${\mathbb E}[T_n^2]$, which then translate into a system of functional equations:
\[
S_2 \equiv \left\{
\begin{array}{ll}
\tilde{S}(z)-z=\tilde{S}(z)^2+2\tilde{R}(z)^2+z^2C(z)^2\\
\tilde{V}(z)=2\tilde{V}(z) \, \tilde{R}(z)+2zC(z) \, T(z)+\tilde{S}(z)-\tilde{R}(z)\\
U(z)=2zC(z) \, U(z)+2T(z)^2+2V(z)-S(z).\\
\end{array}
\right.
\]
Focusing on the first equation in $S_2$, we observe that eq.~(\ref{sec-rn}) gives for $n \geq 2$
$$C_{n-1} \mathbb{E}[\tilde{R}_n^2] = \sum_{j=1}^{n-1} C_{j-1} C_{n-1-j} \Big(\mathbb{E}[\tilde{R}_j^2] \, \mathbb{E}[\tilde{R}_{n-j}^2] + 2 \mathbb{E}[\tilde{R}_j] \, \mathbb{E}[\tilde{R}_{n-j}] + 1\Big), $$
which multiplied by $z^n$ can be rewritten as
\begin{eqnarray*}
C_{n-1} \mathbb{E}[\tilde{R}_n^2] z^n &=& \sum_{j=1}^{n-1} C_{j-1} \mathbb{E}[\tilde{R}_j^2] z^j \cdot C_{n-1-j} \mathbb{E}[\tilde{R}_{n-j}^2] z^{n-j} + 2 \sum_{j=1}^{n-1} C_{j-1} \mathbb{E}[\tilde{R}_j] z^j \cdot C_{n-1-j} \mathbb{E}[\tilde{R}_{n-j}] z^{n-j} \\
&& + z^2 \sum_{j=1}^{n-1} C_{j-1} z^{j-1} \cdot C_{n-1-j} z^{n-1-j}.
\end{eqnarray*}
Summing over $n \geq 2$, we obtain
\begin{eqnarray*}
\tilde{S}(z) &=& z + \sum_{n=2}^\infty C_{n-1} \mathbb{E}[\tilde{R}_n^2] z^n = z + \sum_{n=2}^\infty \sum_{j=1}^{n-1} C_{j-1} \mathbb{E}[\tilde{R}_j^2] z^j \cdot C_{n-1-j} \mathbb{E}[\tilde{R}_{n-j}^2] z^{n-j} \\
&& + 2 \sum_{n=2}^\infty \sum_{j=1}^{n-1} C_{j-1} \mathbb{E}[\tilde{R}_j] z^j \cdot C_{n-1-j} \mathbb{E}[\tilde{R}_{n-j}] z^{n-j}
 + z^2 \sum_{n=2}^\infty \sum_{j=1}^{n-1} C_{j-1} z^{j-1} \cdot C_{n-1-j} z^{n-1-j} \\
 &=& z + \tilde{S}(z)^2+2\tilde{R}(z)^2+z^2C(z)^2.
\end{eqnarray*}
Similarly, the second equation of $S_2$ follows from eq.~(\ref{sec-rnn}), and the third equation from eq.~(\ref{sec-rnnn}).

To solve system $S_2$, we first find $R(z)$ from the first equation of system $S_1$ from the proof of Proposition \ref{meantotuni}. We then find $\tilde{R}(z)$ by using eq.~(\ref{erretilde}). From the first equation of system $S_2$, we can obtain $\tilde{S}(z)$ and also $S(z)$ from eq.~(\ref{esse}). $\tilde{S}(z)$ is then used together with $T(z)$ (eq.~\ref{eqtt}) for calculating $ \tilde{V}(z)$ from the second equation of $S_2$. Once we have a formula for $ \tilde{V}(z)$, we obtain $V(z)$ from eq.~(\ref{vu}), and we finally compute $U(z)$ from the third equation of $S_2$. Writing $r \equiv \sqrt{1-4z}$, we find
\begin{small}
\begin{eqnarray}\label{evv}
%
%
V(z)&=&\frac{-\sqrt{2 r-1}+r \left(-r+\sqrt{2 r-1}-\sqrt{-2r+4 \sqrt{2 r-1}-1}+3\right)-1}{2 r \sqrt{2 r-1}},\\\label{euu}
U(z)&=&\frac{1}{2} \left(-\frac{1}{r^3}+\frac{-\frac{2 \sqrt{-2r+4 \sqrt{2 r-1}-1}}{\sqrt{2 r-1}}+\sqrt{-2 r+4 \sqrt{2 r-1}-1}+\frac{4}{\sqrt{2 r-1}}+3}{r}-\frac{6}{\sqrt{2 r-1}}+1\right).
\end{eqnarray}
\end{small}

The dominant singularity of the generating functions $V(z)$ and $U(z)$ is at $\alpha \equiv 7(8\sqrt{2}-11)/16$, which is the dominant singularity of the square root $\sqrt{-2r+4 \sqrt{2 r-1}-1}$ appearing in eqs.~(\ref{evv}) and (\ref{euu}). We obtain the expansion of $V(z)$ and $U(z)$ at their dominant singularity $\alpha$ by plugging the expansion
\begin{equation*}
\sqrt{-2r+4 \sqrt{2 r-1}-1} \stackrel{z \rightarrow \alpha}{\sim} \sqrt{\frac{7(11-\sqrt{2})}{34}} \cdot \sqrt{1-\frac{16 z}{7(8 \sqrt{2}-11)}}
\end{equation*}
in eqs.~(\ref{evv}) and (\ref{euu}), while setting $z = \alpha$ elsewhere. Algebraic manipulations then lead to
\begin{eqnarray*}
V(z) &\stackrel{z \rightarrow \alpha}{\sim}& k_1 - \frac{1}{2}\left(1+\frac{\sqrt{2}}{2}\right)\sqrt{\frac{7(11-\sqrt{2})}{34}} \cdot \sqrt{1-\frac{16 z}{7(8 \sqrt{2}-11)}},\\
U(z) &\stackrel{z \rightarrow \alpha}{\sim}& k_1' - \frac{1}{17}(15+11\sqrt{2})\sqrt{\frac{7(11-\sqrt{2})}{34}} \cdot \sqrt{1-\frac{16 z}{7(8 \sqrt{2}-11)}} ,
\end{eqnarray*}
for certain constants $k_1$ and $k_1'$. Eq.~(\ref{pillo2}), together with the asymptotic expansion in eq.~(\ref{cata}), finally yields
\begin{eqnarray} \label{ett1}
{\mathbb E}[T_n R_n] &=& \frac{[z^n] V(z)}{C_{n-1}} \sim \left(1+\frac{\sqrt{2}}{2}\right)\sqrt{\frac{7(11-\sqrt{2})}{34}}\left[\frac{4}{7 (8\sqrt{2} - 11) }\right]^n, \\ \label{ett2}
{\mathbb E}[T_n^2] &=& \frac{[z^n] U(z)}{C_{n-1}} \sim \frac{2}{17}(15+11\sqrt{2})\sqrt{\frac{7(11-\sqrt{2})}{34}}\left[\frac{4}{7 (8\sqrt{2} - 11) }\right]^n.
\end{eqnarray}

By using these calculations, we obtain the following result.
\begin{prop}\label{prop:varunif}
The variance of the total number $T_n$ of ancestral configurations in an ordered unlabeled topology of size $n$ selected uniformly at random satisfies the asymptotic relation ${\mathbb V}[T_n] \sim {\mathbb E}[T_n^2]$, where ${\mathbb E}[T_n^2]$ grows as in eq.~(\ref{ett2}). For increasing values of $n$, the correlation coefficient $\rho[T_n,R_n]$ between the total number $T_n$ of ancestral configurations and the number $R_n$ of root configurations converges to a constant $\rho[T_n,R_n] \rightarrow 0.9004\ldots.$
\end{prop}
\noindent \emph{Proof.}
First, for the variance we have ${\mathbb V}[T_n] = \mathbb{E}[T_n^2] - \mathbb{E}[T_n]^2 \sim{\mathbb E}[T_n^2]$. Indeed, from Proposition \ref{meantotuni}, $\mathbb{E}[T_n]^2 \bowtie [ ({4}/{3})^2 ]^n = ({16}/{9})^n$, whereas from eq.~(\ref{ett2}), $\mathbb{E}[T_n^2] \bowtie \big[{4}/[7 (8\sqrt{2} - 11)]\big]^n$, and ${4}/{[7 (8\sqrt{2} - 11)]} > \frac{16}{9}$.

Second, the covariance ${\rm Cov}[T_n,R_n]$ grows like ${\rm Cov}[T_n,R_n] = \mathbb{E}[T_n R_n] - \mathbb{E}[T_n] \, \mathbb{E}[R_n] \sim {\mathbb E}[T_nR_n]$. Indeed, by Proposition \ref{meantotuni} and eq.~(\ref{meandis-propyu}) for the uniform model, we have $\mathbb{E}[T_n] \, \mathbb{E}[R_n] \bowtie ({4}/{3})^n \cdot (4/3)^n = ({16}/{9})^n$; from eq.~(\ref{ett1}), we have $\mathbb{E}[T_n R_n] \bowtie \big[ {4}/[7 (8\sqrt{2} - 11)] \big]^n$, with ${4}/[7 (8\sqrt{2} - 11)] > \frac{16}{9}$. Hence, for the correlation coefficient,
\[
\rho[T_n,R_n] = \frac{{\rm Cov}[T_n,R_n]}{\sqrt{{\mathbb V}[T_n]} \sqrt{{\mathbb V}[R_n]}} \sim \frac{{\displaystyle 1+\frac{\sqrt{2}}{2}}}{{\displaystyle \sqrt{\frac{2}{17}(15+11\sqrt{2})}}} \approx 0.9004,
\]
where we used the asymptotic formula for ${\mathbb V}[R_n]$ given for the uniform model in eq.~(\ref{vardis-propyuv}) as well as the asymptotics in eqs.~(\ref{ett1}) and (\ref{ett2}) for ${\rm Cov}[T_n,R_n]$ and ${\mathbb V}[T_n]$, respectively.
$\Box$

\vskip .2cm
To conclude this section, we show that the total number $T_n$ of ancestral configurations of an ordered unlabeled topology of size $n$ selected uniformly at random has an asymptotic lognormal distribution. From Lemma \ref{lemcorr} and Section \ref{distrib}, we know that the logarithm of the number $R_n$ of root configurations in an ordered unlabeled topology of size $n$ selected uniformly at random converges asymptotically to a standard normal distribution, that is,
\begin{equation}
\label{eq:normal}
\frac{\log R_n-{\mathbb E}[\log R_n]}{ \sqrt{ {\mathbb V}[\log R_n]} } \stackrel{d}{\longrightarrow}N(0,1),
\end{equation}
with ${\mathbb E}[\log R_n] \sim 0.272 \cdot n $ and ${\mathbb V}[\log R_n] \sim 0.034 \cdot n$ (where the constants are approximate).

\textcolor{black}{From eq.~(\ref{bau}), the variables $R_n$ and $T_n$ measured over the same random ordered unlabeled topology of size $n$ satisfy $R_n\leq T_n \leq (2n-1)R_n$. This inequality gives $\log T_n = \log R_n + \epsilon_n$, where the random variable $\epsilon_n$ has values in $[0, \log(2n-1)]$. Thus, we have
\begin{eqnarray}\label{logN1}
{\mathbb E}[\log T_n] &=& {\mathbb E}[\log R_n]+ {\mathbb E}[\epsilon_n] = {\mathbb E}[\log R_n]+{\mathcal O}(\log n) \sim {\mathbb E}[\log R_n] \\
\label{logN2}
{\mathbb V}[\log T_n] &=& {\mathbb V}[\log R_n] + {\mathbb V}[\epsilon_n] + 2 {\rm Cov}[\log R_n, \epsilon_n]
\sim {\mathbb V}[\log R_n],
\end{eqnarray}
where we use ${\mathbb V}[\epsilon_n] \leq [\log(2n-1)]^2/4$ from Popoviciu's inequality on the maximal variance for a bounded random variable~\cite{BhatiaAndDavis00}, so that the comparison with the linearly increasing ${\mathbb V}[\log R_n]$, gives $\lim_{n \rightarrow \infty} {\mathbb V}[\epsilon_n] / {\mathbb V}[\log R_n]= 0$; we also use ${\rm Cov}[\log R_n, \epsilon_n] \leq \sqrt{{\mathbb V}[\log R_n]} \, \sqrt{{\mathbb V}[\epsilon_n]}$ from the Cauchy-Schwarz inequality.}

\textcolor{black}{Next, we write
\begin{eqnarray}
\label{eq:slutsky}
\frac{\log T_n-{\mathbb E}[\log T_n]}{\sqrt{{\mathbb V}[\log T_n]}} & = & \frac{\log R_n-{\mathbb E}[\log R_n]}{\sqrt{{\mathbb V}[\log R_n]}}\cdot \frac{\sqrt{{\mathbb V}[\log R_n]}}{\sqrt{{\mathbb V}[\log T_n]}} + \frac{\epsilon_n-{\mathbb E}[\epsilon_n]}{\sqrt{{\mathbb V}[\log T_n]}}.
\end{eqnarray}
The expression $({\log R_n-{\mathbb E}[\log R_n]})/{\sqrt{{\mathbb V}[\log R_n]}}$ converges in distribution to a normal random variable with mean 0 and variance 1 (eq.~\ref{eq:normal}). The ratio ${\sqrt{{\mathbb V}[\log R_n]}}/{\sqrt{{\mathbb V}[\log T_n]}}$ is a number sequence that by eq.~(\ref{logN2}) converges to a finite constant, 1. The expression $({\epsilon_n-{\mathbb E}[\epsilon_n]}/{\sqrt{{\mathbb V}[\log T_n]}})$ converges in mean square to 0, as
$$\lim_{n \rightarrow \infty} {\mathbb E}\bigg[ \bigg( \frac{\epsilon_n-{\mathbb E}[\epsilon_n]}{\sqrt{{\mathbb V}[\log T_n]}} - 0 \bigg)^2 \bigg] = \lim_{n \rightarrow \infty} \frac{ {\mathbb V}[\epsilon_n]}{ {\mathbb V}[\log T_n] };$$
the denominator ${\mathbb V}[\log T_n]$ increases linearly with $n$ (Section \ref{distrib}), and again by Popoviciu's inequality, the numerator is bounded above by $[\log(2n-1)]^2/4$, so that $\lim_{n \rightarrow \infty} { {\mathbb V}[\epsilon_n]} / { {\mathbb V}[\log T_n] } = 0$.}

\textcolor{black}{As convergence in mean square implies convergence in probability~(see p.~10 of \cite{Serfling80}), we can apply Slutsky's theorem on perturbation of random variables that converge in distribution by random variables that converge in probability~(see p.~19 of \cite{Serfling80}) to eq.~(\ref{eq:slutsky}). In particular, the convergence in distribution of $({\log R_n-{\mathbb E}[\log R_n]}) \, \times$$\\({\sqrt{{\mathbb V}[\log R_n]}})^{-1}$, trivial convergence in probability of  ${\sqrt{{\mathbb V}[\log R_n]}}/{\sqrt{{\mathbb V}[\log T_n]}}$, and convergence in probability of $({\epsilon_n-{\mathbb E}[\epsilon_n]} )/{\sqrt{{\mathbb V}[\log T_n]}}$ allow us to conclude $({\log T_n-{\mathbb E}[\log T_n]})/{\sqrt{{\mathbb V}[\log T_n]}}$ converges in distribution to a normal random variable with mean 0 and variance 1.}

\begin{figure}
\begin{center}
\includegraphics*[scale=0.40,trim=0 0 0 0]{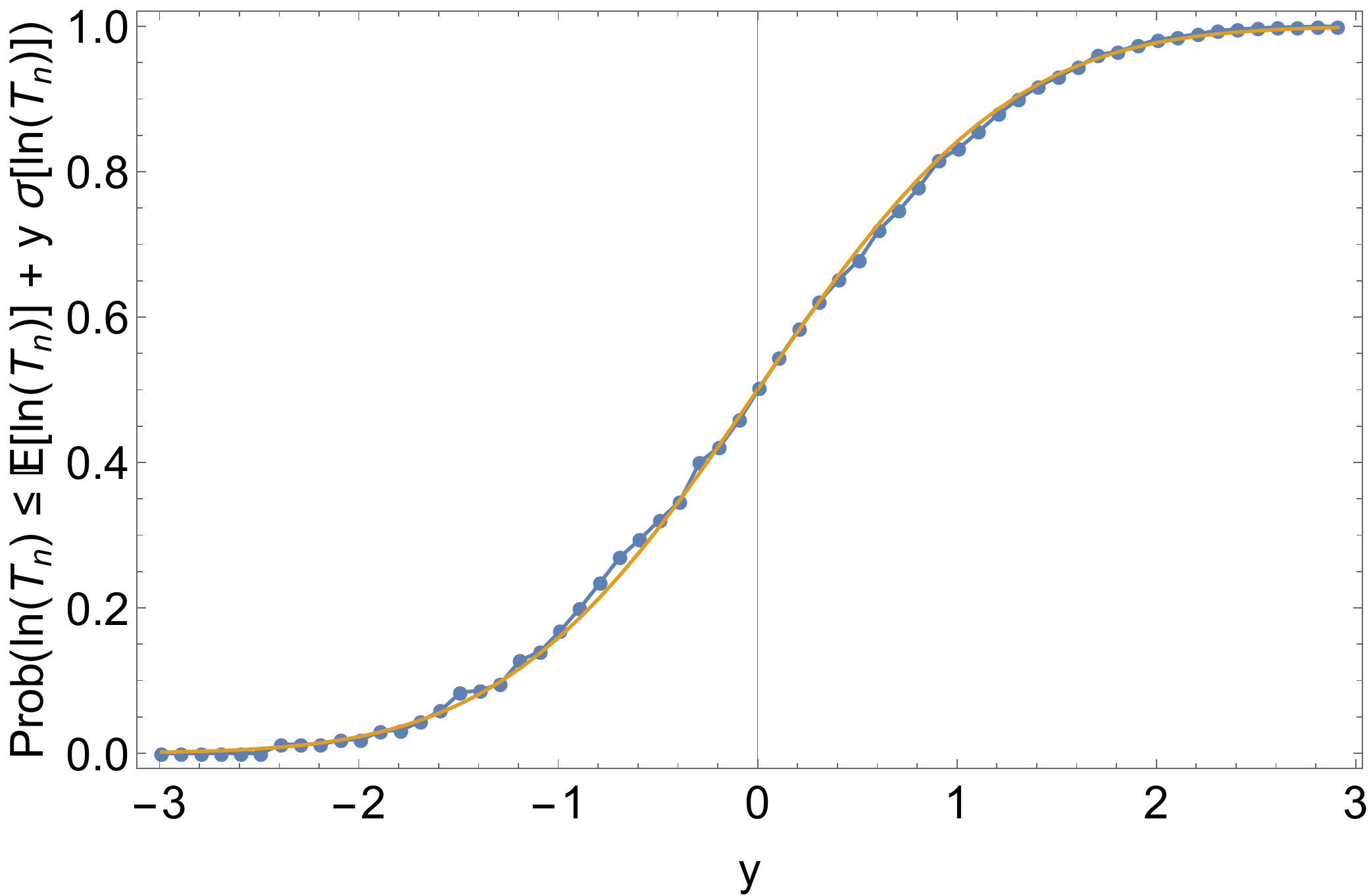}
\end{center}
\vspace{-.7cm}
\caption{{\small Cumulative distribution of the natural logarithm of the total number of configurations for uniformly distributed ordered unlabeled topologies (or uniformly distributed labeled topologies) of size $n = 15$ (dotted line). For each $y \in [-3, 3]$ in steps of $0.1$, the quantity plotted is the probability that an ordered unlabeled topology (or labeled topology) with $n = 15$ chosen uniformly at random has total number of configurations less than or equal to $\exp (\mathbb{E}[\log T_n] + y \sigma[\log T_n])$, where $\mathbb{E}[\log T_n]$ and $\sigma[\log T_n]=\sqrt{{\mathbb V}[\log T_n]}$ are respectively the mean and standard deviation of the logarithm of the total number of configurations for uniformly distributed ordered unlabeled topologies (or labeled topologies) with $n = 15$ taxa. The solid line is the cumulative distribution of a Gaussian random variable with mean 0 and variance 1.
}} \label{normalunif}
\end{figure}

Fig.~\ref{normalunif} shows the cumulative distribution $\mathbb{P}[\log T_n \leq \mathbb{E}[\log T_n] + y \sqrt{{\mathbb V}[\log T_n]}]$ as a function of $y$, when ordered unlabeled topologies of size 15 are selected uniformly at random.  To obtain the distribution, we count total configurations for each unlabeled topology $t$ with 15 leaves, and then count the number of ordered unlabeled topologies having the shape of $t$ (eq.~\ref{outeq}). The figure illustrates the agreement between the exact cumulative distribution of ancestral configurations and the standard normal distribution.

By the equivalence between ordered unlabeled topologies and labeled topologies reported in Lemma \ref{lemcorr}, we can state the main results of this section as follows.
\begin{theo}\label{totttuni}
For a labeled topology of size $n$ selected at random under the uniform distribution, the mean and the variance of the total number $T_n$ of ancestral configurations grow asymptotically like
\begin{eqnarray}
\mathbb{E}[T_n] &\sim & \sqrt{6} \left(\frac{4}{3}\right)^n,\\
\mathbb{V}[T_n] &\sim & \frac{2}{17}(15+11\sqrt{2})\sqrt{\frac{7(11-\sqrt{2})}{34}}\left[\frac{4}{7 (8\sqrt{2} - 11) }\right]^n.
\end{eqnarray}
Furthermore, the logarithm of the total number of ancestral configurations in a labeled topology of size $n$ selected uniformly at random, rescaled as $ ({\log T_n-{\mathbb E}[\log T_n]})/{\sqrt{{\rm Var}[\log T_n]}}$, converges to a standard normal distribution, where ${\mathbb E}[\log T_n] \sim \mu n $ and ${\mathbb V}[\log T_n] \sim \sigma^2 n$, $(\mu, \sigma^2) \approx (0.272, 0.034)$.
\end{theo}

\subsection{Uniform ordered unlabeled histories and Yule labeled topologies}

By Lemma \ref{lemcorr}, the distribution of the total number of ancestral configurations over random labeled topologies of size $n$ selected under the Yule probability model is the distribution of the total number of configurations over uniformly distributed ordered unlabeled histories of size $n$. We exploit  this equivalence for this section, now denoting by $R_n$ and $T_n$, respectively, the number of root ancestral configurations and the total number of ancestral configurations in a random ordered unlabeled history of size $n$ under the uniform distribution.
\begin{lemm}\label{rectotyu}
Consider the random variables $R_n$ and $T_n$. We have $R_1=T_1=0$, and for $n\geq 2$,
\begin{eqnarray}
\label{eq:Rn}
R_n    & \stackrel{d}{=} & R_{I_n} \, R^{*}_{n-I_n} + R_{I_n} + R^{*}_{n-I_n} + 1, \\
\label{eq:Tn}
T_n    & \stackrel{d}{=} & T_{I_n}+T_{n-I_n}^{*}+R_n, \\\nonumber
T_nR_n &\overset{d}{=}   &
T_{I_n}R_{I_n}R_{n-I_n}^* +
T_{I_n}R_{I_n}+
T_{I_n}R_{n-I_n}^* +
T_{I_n} \\\label{eq:TnRn} && + \,
T_{n-I_n}^*R_{I_n}R_{n-I_n}^* +
T_{n-I_n}^*R_{I_n}+
T_{n-I_n}^*R_{n-I_n}^* +
T_{n-I_n}^* +
(R_n)^2,\\
\label{eq:Tn2}
(T_n)^2 &\overset{d}{=} & (T_{I_n})^2+(T_{I_{n-I_n}}^*)^2+2T_{I_n}T_{n-I_n}^*+2T_nR_n-(R_n)^2,
\end{eqnarray}
where $I_n$ is a uniformly distributed variable over the interval $[1,n-1]$, $R^{*}_{j}$ and $T^{*}_{j}$ are independent copies of $R_j$ and $T_j$, respectively, for every $j \in [1,n-1]$, and both $R_j$ and $R_j^{*}$ as well as $T_j,$ and $T_j^{*} $ are independent of $I_j$ for $j \in [1,n-1]$.
\end{lemm}
\noindent \emph{Proof.} Eqs.~(\ref{eq:Rn}) and (\ref{eq:Tn}) are from Lemma \ref{rectot}. By expanding $T_nR_n \overset{d}{=} (T_{I_n}+T_{n-I_n}^*)(R_{I_n}+1)(R_{n-I_n}^*+1) +(R_n)^2$, we have eq.~(\ref{eq:TnRn}). Finally, eq.~(\ref{eq:Tn2}) is obtained by squaring eq.~(\ref{eq:Tn}); it also copies eq.~(\ref{sec-rnnn}) from Lemma \ref{4rec}. $\Box$

\vskip .2cm
The distributional recurrences in Lemma \ref{rectotyu} can be used to determine recurrences for the expectations ${\mathbb E}[R_n]$, ${\mathbb E}[T_n]$, ${\mathbb E}[T_nR_n]$, and ${\mathbb E}[T_n^2]$. For $n \geq 1$, we can write
\begin{eqnarray}\label{recurr0}
(n-1) \mathbb{E}[R_n]        &=&  \bigg( \sum_{j=1}^{n-1}\mathbb{E}[R_j] \, \mathbb{E}[R_{n-j}] \bigg) + 2 \bigg( \sum_{j=1}^{n-1} \mathbb{E}[R_j] \bigg) + (n-1), \\\label{recurr1}
(n-1) \mathbb{E}[T_n]        &=& 2 \bigg( \sum_{j=1}^{n-1}\mathbb{E}[T_j] \bigg) +(n-1) \mathbb{E}[R_n],\\\nonumber
(n-1) \mathbb{E}[T_n R_n] &=& 2
\bigg( \sum_{j=1}^{n-1}\mathbb{E}[T_jR_j] \, \mathbb{E}[R_{n-j}] \bigg) +2 \bigg( \sum_{j=1}^{n-1}\mathbb{E}[T_jR_j] \bigg) + 2 \bigg( \sum_{j=1}^{n-1}\mathbb{E}[T_j] \, \mathbb{E}[R_{n-j}] \bigg) \\\label{recurr2}
&&+2 \bigg( \sum_{j=1}^{n-1}\mathbb{E}[T_j] \bigg) +(n-1)\mathbb{E}[R_n^2],\\\label{zizzo}
(n-1)\mathbb{E}[T_n^2]&=&2 \bigg( \sum_{j=1}^{n-1}\mathbb{E}[T_j^2] \bigg) +2 \bigg( \sum_{j=1}^{n-1}\mathbb{E}[T_j] \, \mathbb{E}[T_{n-j}] \bigg) +2(n-1)\mathbb{E}[T_nR_n]-(n-1)\mathbb{E}[R_n^2].
\end{eqnarray}
Define the following generating functions
\begin{eqnarray}
\label{eq:Rz}
R(z) & \equiv & \sum_{n=1}^\infty \mathbb{E}[R_n]z^n, \\
\label{eq:Tz}
T(z) & \equiv & \sum_{n=1}^\infty \mathbb{E}[T_n]z^n, \\
\label{eq:Sz}
S(z) & \equiv & \sum_{n=1}^\infty \mathbb{E}[R_n^2]z^n, \\
\label{eq:Vz}
V(z) & \equiv & \sum_{n=1}^\infty \mathbb{E}[T_n R_n]z^n, \\
\label{eq:Uz}
U(z) & \equiv & \sum_{n=1}^\infty \mathbb{E}[T_n^2]z^n.
\end{eqnarray}

The recurrences in eqs.~(\ref{recurr1}), (\ref{recurr2}), and (\ref{zizzo}) translate into a system of differential equations:
\[
S_3 \equiv \left\{
\begin{array}{ll}
T'(z)-\frac{z+1}{z-z^2}T(z) = R'(z)-\frac{R(z)}{z}\\
V'(z)-\left(\frac{2R(z)}{z}+\frac{z+1}{z-z^2}\right)V(z) = \frac{2T(z) \, R(z)+\frac{2z}{1-z}T(z)+zS'(z)-S(z)}{z}\\
U'(z)-\frac{z+1}{z-z^2}U(z) = \frac{2T(z)^2+2zV'(z)-2V(z)-zS'(z)+S(z)}{z},\\
\end{array}
\right.
\]
The derivatives $R'(z)$, $T'(z)$, $S'(z)$, $V'(z)$, and $U'(z)$ appear in $S_3$ due to the factor $n-1$ in eqs.~(\ref{recurr1}), (\ref{recurr2}), and (\ref{zizzo}). We derive the third equation in $S_3$ as an example. First, multiplying both sides of eq.~(\ref{zizzo}) by $z^n$ yields
\begin{eqnarray}\nonumber
 z n \mathbb{E}[T_n^2] z^{n-1} - \mathbb{E}[T_n^2] z^{n} &=& 2 \sum_{j=1}^{n-1}\mathbb{E}[T_j^2] z^j \cdot z^{n-j} + 2 \sum_{j=1}^{n-1}\mathbb{E}[T_j] z^j \cdot \mathbb{E}[T_{n-j}] z^{n-j} + 2 z n \mathbb{E}[T_n R_n] z^{n-1} \\\nonumber
 && - 2 \mathbb{E}[T_n R_n] z^n - z n \mathbb{E}[R_n^2] z^{n-1} + \mathbb{E}[R_n^2] z^n.
\end{eqnarray}
Summing over $n \geq 1$, we obtain
\begin{eqnarray}\nonumber
 z \sum_{n=1}^\infty n \mathbb{E}[T_n^2] z^{n-1} - \sum_{n=1}^\infty \mathbb{E}[T_n^2] z^{n} &=& 2 \sum_{n=1}^\infty \sum_{j=1}^{n-1}\mathbb{E}[T_j^2] z^j \cdot z^{n-j} + 2 \sum_{n=1}^\infty \sum_{j=1}^{n-1}\mathbb{E}[T_j] z^j \cdot \mathbb{E}[T_{n-j}] z^{n-j} \\\nonumber && + 2 z \sum_{n=1}^\infty n \mathbb{E}[T_n R_n] z^{n-1}
 - 2 \sum_{n=1}^\infty \mathbb{E}[T_nR_n] z^n - z \sum_{n=1}^\infty n \mathbb{E}[R_n^2] z^{n-1} + \sum_{n=1}^\infty \mathbb{E}[R_n^2] z^n.
\end{eqnarray}
To complete the derivation, we note that this equation can be rewritten:
\begin{equation}\nonumber
 z U'(z) - U(z) = 2 U(z) \, \left( \frac{z}{1-z} \right) + 2 T(z)^2 + 2 z V'(z) - 2 V(z) - z S'(z) + S(z),
\end{equation}
as $U'(z) = \left( \sum_{n=1}^\infty \mathbb{E}[T_n^2]z^n \right)' = \sum_{n=1}^\infty n \mathbb{E}[T_n^2] z^{n-1}$, $V'(z) = \left( \sum_{n=1}^\infty \mathbb{E}[T_n R_n]z^n \right)' = \sum_{n=1}^\infty n \mathbb{E}[T_n R_n] z^{n-1}$, $S'(z) = \left( \sum_{n=1}^\infty \mathbb{E}[R_n^2]z^n \right)' = \sum_{n=1}^\infty n \mathbb{E}[R_n^2] z^{n-1}$, and $\frac{z}{1-z} = \sum_{n=1}^\infty z^n$.

We also observe that the generating functions $R(z)$ and $S(z)$ (eqs.~\ref{eq:Rz} and \ref{eq:Sz}) were studied in the analysis of root configurations under the Yule model for labeled topologies in Section 5 of \cite{DisantoEtAl20}. In particular, eq.~(39) in the proof of Proposition 5.3 of \cite{DisantoEtAl20} found that $R(z)$---there denoted by $E(z)$---has explicit form
\begin{equation}\label{formE}
R(z) = \frac{2 z \sin\left( \frac{\sqrt{3}}{2} \log(1-z) \right)}{(z-1)\left[ \sqrt{3} \cos\left( \frac{\sqrt{3}}{2} \log(1-z) \right) + \sin\left( \frac{\sqrt{3}}{2} \log(1-z) \right) \right]}.
\end{equation}
The dominant singularity is $\alpha_1 \equiv 1-e^{-2\pi\sqrt{3}/9}$, and the singular expansion at the dominant singularity is
\begin{equation}\label{singR}
R(z) \overset{z \rightarrow \alpha_1}{\sim} \frac{1}{1- \frac{z}{\alpha_1}}.
\end{equation}

The generating function $S(z)$ was found in Section 5.3 of \cite{DisantoEtAl20} to have singular expansion
\begin{equation}\label{asis}
S(z) \overset{z \rightarrow \alpha_2}{\sim} \frac{1}{1 - \frac{z}{\alpha_2}},
\end{equation}
where the dominant singularity $\alpha_2 \equiv 0.4889986317\ldots$ was approximated in the Appendix. By singularity analysis (eq.~\ref{pillo1}), the expansions in eqs.~(\ref{singR}) and (\ref{asis}) yield the asymptotic relations in eqs.~(\ref{meandis-propyu}) and (\ref{vardis-propyuv}):
\begin{equation}\label{ripet}
\mathbb{E}[R_n] \sim \alpha_1^{-n} \quad \text{and} \quad \mathbb{E}[R_n^2] \sim \alpha_2^{-n}.
\end{equation}
Note indeed, that the asymptotic constant $2.0449954971\ldots$ appearing in eq.~(\ref{vardis-propyuv}) is obtained as $\alpha_2^{-1}$.

We now observe that eq.~(\ref{singR}) and the first equation of $S_3$ yield the asymptotic growth of the mean number $\mathbb{E}[T_n]$ of ancestral configurations in an ordered unlabeled history of size $n$ selected uniformly at random.
\begin{prop}\label{proppa}
The mean total number of ancestral configurations in an ordered unlabeled history of size $n$ selected uniformly at random satisfies the asymptotic relation $\mathbb{E}[T_n] \sim \mathbb{E}[R_n] \sim \alpha_1^{-n}= [{1}/({1-e^{-2\pi\sqrt{3}/9}})]^{n}.$
\end{prop}
\noindent \emph{Proof.} We start by rewriting the first equation of $S_3$ as
\begin{equation}
\label{eq:TprimeM}
T'(z) \, M(z)-\frac{z+1}{z-z^2}M(z) \, T(z)=\left[R'(z)-\frac{R(z)}{z}\right]M(z),
\end{equation}
where $M(z)\equiv {(z-1)^2}/{z}$ is the integrating factor.

Since $M'(z) = -\frac{z+1}{z-z^2}M(z)$, the left-hand side of eq.~(\ref{eq:TprimeM}) can be rewritten $[T(z) \, M(z)]'$, yielding
\begin{equation}
\label{eq:TMprime}
\left[T(z)\frac{(z-1)^2}{z}\right]'=\left[R'(z)-\frac{R(z)}{z}\right]\frac{(z-1)^2}{z}.
\end{equation}
Because $T_1=0$, the expansion of $T(z)$ starts with a non-zero quadratic term. Hence, we have
$$\left[T(z)\frac{(z-1)^2}{z}\right]_{z=0} = 0,$$
and the differential equation in eq.~(\ref{eq:TMprime}) thus gives $T(z)\frac{(z-1)^2}{z} = \int_{0}^{z} [R'(t)- {R(t)}/{t}]\frac{(t-1)^2}{t} \, dt$, that is,
\[
T(z)=\frac{z}{(z-1)^2}\int_{0}^{z}\left[R'(t)-\frac{R(t)}{t}\right]\frac{(t-1)^2}{t} \, dt.
\]

To obtain the singular expansion of $T(z)$, we must analyze functions $R'(t)$, $[R'(t)-{R(t)}/{t}]\frac{(t-1)^2}{t}$, and $\int_{0}^{z}[R'(t)-{R(t)}/{t}]\frac{(t-1)^2}{t} \, dt$ at their dominant singularity. Because $\alpha_1 = 1-e^{-2 \pi \sqrt{3}/9}$ is the dominant singularity of $R(t)$ and $R(t)\overset{t\rightarrow \alpha_1}{\sim} {1}/({1-\frac{t}{\alpha_1}})$ (eq.~\ref{singR}), from eq.~(\ref{diffe}), $R'(t)$ has dominant singularity at $t = \alpha_1$. Its singular expansion is
$$R'(t)\overset{t\rightarrow \alpha_1}{\sim}\frac{1}{\alpha_1(1-\frac{t}{\alpha_1})^2},$$
obtained by differentiating the expansion of $R(t)$. It follows that $\alpha_1$ is also the dominant singularity of the function $[R'(t)-{R(t)}/{t}]\frac{(t-1)^2}{t}$, whose singular expansion follows
\begin{equation}\label{kiju}
\left[R'(t)-\frac{R(t)}{t}\right]\frac{(t-1)^2}{t} \overset{t\rightarrow \alpha_1}{\sim} \left[ \frac{1}{\alpha_1(1-\frac{t}{\alpha_1})^2} \right]\frac{(\alpha_1-1)^2}{\alpha_1}.
\end{equation}
Finally, by eq.~(\ref{inte}) $\int_{0}^{z}[R'(t)-{R(t)}/{t}] \frac{(t-1)^2}{t} \, dt$ can be expanded at its dominant singularity $\alpha_1$ by integrating the singular expansion of the integrand function (eq.~\ref{kiju}). Consequently, the expansion of $T(z)$ at its dominant singularity $\alpha_1$ satisfies
\begin{equation*}
T(z) \overset{z\rightarrow \alpha_1}{\sim}\frac{\alpha_1}{(\alpha_1-1)^2}\int_{0}^{z}\frac{1}{\alpha_1(1-\frac{t}{\alpha_1})^2}\frac{(\alpha_1-1)^2}{\alpha_1} \, dt = \frac{1}{1-\frac{z}{\alpha_1}} -1
\overset{z\rightarrow \alpha_1}{\sim}
\frac{1}{1-\frac{z}{\alpha_1}}.
\end{equation*}
By eqs.~(\ref{pillo1}) and (\ref{ripet}), we conclude
\[
\mathbb{E}[T_n]=[z^n]T(z) \sim \alpha_1^{-n} \sim \mathbb{E}[R_n].
 \quad \Box \]

\begin{figure}
\begin{center}
\includegraphics*[scale=0.40,trim=0 0 0 0]{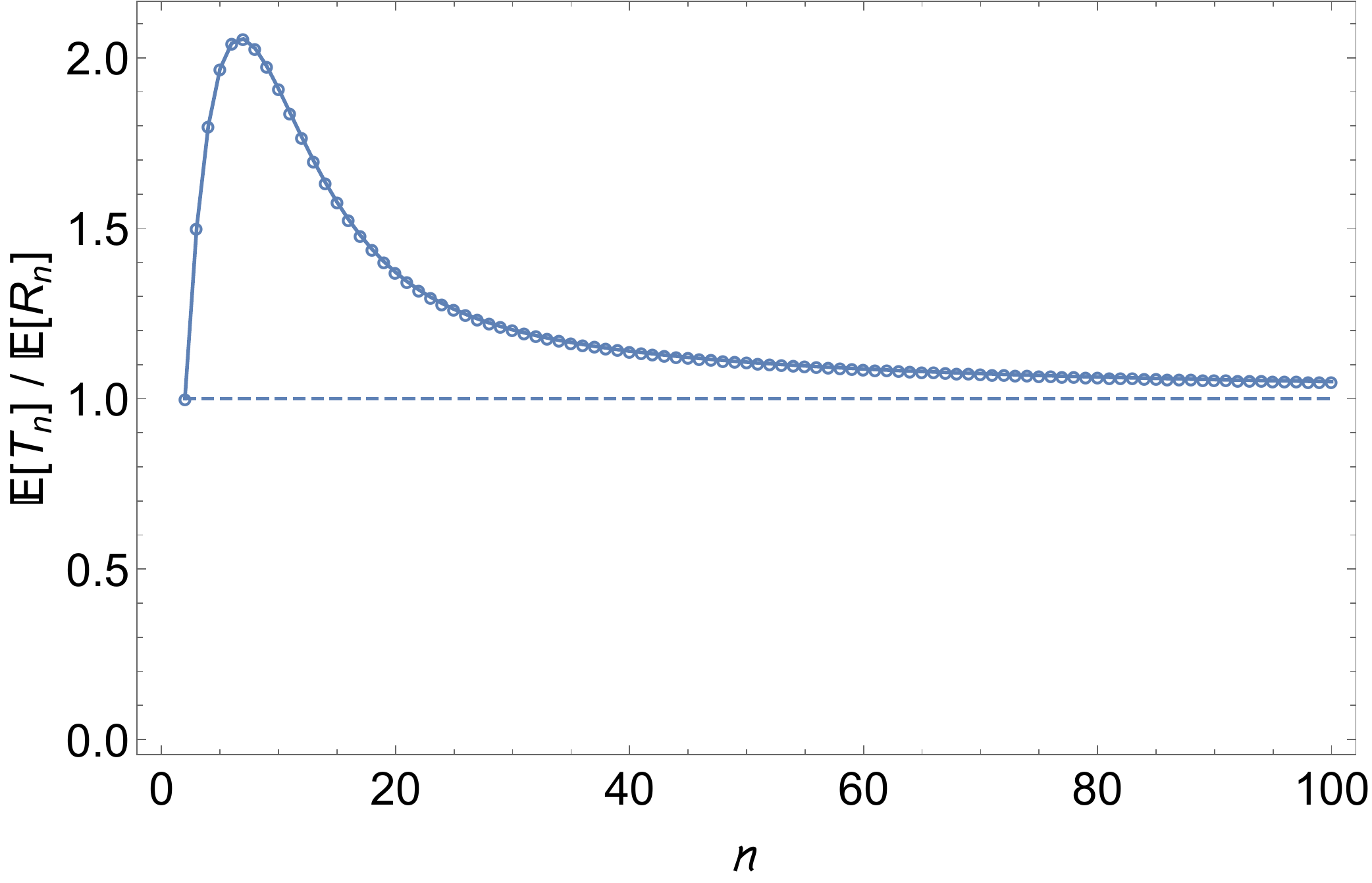}
\end{center}
\vspace{-.7cm}
\caption{{\small Ratio of the mean total number $\mathbb{E}[T_n]$ of ancestral configurations and mean number $\mathbb{E}[R_n]$ of root configurations for uniformly distributed ordered unlabeled histories (or Yule-distributed labeled topologies) of size $2 \leq n \leq 100$. Values of $\mathbb{E}[R_n]$ and $\mathbb{E}[T_n]$ are computed from the recurrences in eqs.~(\ref{recurr0}) and (\ref{recurr1}), respectively. }
} \label{ratioyule}
\end{figure}

\vskip .2cm
In Fig.~\ref{ratioyule}, we show a numerical plot of the ratio $\mathbb{E}[T_n] / \mathbb{E}[R_n]$ as a function of $n$. Following the proposition, as $n$ increases, the numerical ratio approaches 1.

We now study the variance ${\mathbb V}[T_n]$ of the total number of ancestral configurations and the correlation coefficient $\rho [T_n, R_n]$ between the numbers of total and root configurations in uniformly distributed ordered unlabeled histories of fixed size. By again using properties of singular expansions under differentiation and integration, we determine the asymptotic growth of the expectations ${\mathbb E}[T_nR_n]$ and ${\mathbb E}[T_n^2]$. We start with ${\mathbb E}[T_nR_n]$. We abbreviate a term in the second equation of $S_3$ by $P(z)$:
\[
P(z)\equiv \frac{2T(z)R(z)+\frac{2z}{1-z}T(z)+zS'(z)-S(z)}{z}.
\]
The second equation of $S_3$ becomes
$V'(z)-V(z)\,[\frac{2R(z)}{z}+\frac{z+1}{z-z^2}] = P(z).$
We introduce integration factor $M(z)$,
\[
M(z)\equiv \left[ e^{\int_0^z -\frac{2R(t)}{t} \, dt} \right] \cdot
\frac{(z-1)^2}{z},
\]
such that $M'(z) = - [\frac{2 R(z)}{z} + \frac{z+1}{z-z^2} ] M(z)$. We find $[V(z) \, M(z)]' = P(z) \, M(z)$, and thus
\[
V(z)=\frac{1}{M(z)}\int_{0}^{z}P(t) \, M(t) \, dt.
\]

To determine the singular expansion of $V(z)$ at its dominant singularity, we observe that $P(t) \, M(t)$ is a function of $T(t)$, $R(t)$, $S(t)$, and $S'(t)$. As demonstrated in Proposition \ref{proppa}, $T(t)$ and $R(t)$ have the same dominant singularity $\alpha_1 \approx 0.702$ (eq.~\ref{singR}), a value larger than the dominant singularity $\alpha_2 \approx 0.489$ of $S(t)$ and $S'(t)$ (eq.~\ref{asis}). Hence, $R(t)$ and $T(t)$ are analytic functions in a neighborhood of $0$, say $|t| \leq \frac{1}{2}$, that contains $\alpha_2$. As a consequence, we can obtain the singular expansion of $P(t) \, M(t)$ at its dominant singularity $\alpha_2$ by replacing $S(t)$ and $S'(t)$ with their expansions $S(t) \overset{t\rightarrow \alpha_2}{\sim} {1}/({1-\frac{t}{\alpha_2}})$ (eq.~\ref{asis}) and $S'(t) \overset{t\rightarrow \alpha_2}{\sim} [ {1}/({1-\frac{t}{\alpha_2}})]' = {1}/[{\alpha_2 (1-\frac{t}{\alpha_2} )^2}]$, while substituting $t=\alpha_2$ elsewhere. We find
$$P(t) \, M(t)\overset{t\rightarrow \alpha_2}{\sim} \left[ \frac{2T(\alpha_2) \, R(\alpha_2)+\frac{2\alpha_2}{1-\alpha_2}T(\alpha_2)}{\alpha_2} + \frac{1}{\alpha_2(1-\frac{t}{\alpha_2})^2} \right] M(\alpha_2) \overset{t\rightarrow \alpha_2}{\sim}\frac{M(\alpha_2)}{\alpha_2(1-\frac{t}{\alpha_2})^2}.$$
The singular expansion under integration (eq.~\ref{inte}) thus gives
$\int_{0}^{z}P(t)\, M(t) \, dt \overset{z\rightarrow \alpha_2}{\sim} \int_{0}^{z} {M(\alpha_2)}/[{\alpha_2(1-\frac{t}{\alpha_2})^2}] \, dt \overset{z\rightarrow \alpha_2}{\sim} {M(\alpha_2)}/({1-\frac{z}{\alpha_2}})$, from which the singular expansion of $V(z)$ at its dominant singularity $\alpha_2$ is
\[
V(z)\overset{z\rightarrow \alpha_2}{\sim}\frac{1}{M(\alpha_2)} \frac{M(\alpha_2)}{1-\frac{z}{\alpha_2}} = \frac{1}{1-\frac{z}{\alpha_2}}.
\]
Hence, by applying eq.~(\ref{pillo1}) together with eq.~(\ref{ripet}), we have
\begin{equation} \label{x}
\mathbb{E}[T_nR_n] = [z^n] V(z) \sim \alpha_2^{-n} \sim \mathbb{E}[R_n^2].
\end{equation}

We follow the same approach to determine the asymptotic growth of ${\mathbb E}[T_n^2]$. Multiplying both sides of the third equation of $S_3$ by the integrating factor $M(z) \equiv {(z-1)^2}/{z}$ used in the proof of Proposition~\ref{proppa}, we find
\[
U(z)=\frac{z}{(z-1)^2}\int_{0}^{z}\left[\frac{2T(t)^2}{t}+2\left(V'(t)-\frac{V(t)}{t}\right)-\left(S'(t)-\frac{S(t)}{t}\right)\right]\frac{(t-1)^2}{t} \, dt.
\]
We abbreviate
\[
G(t)\equiv\frac{2T(t)^2}{t}+2\left(V'(t)-\frac{V(t)}{t}\right)-\left(S'(t)-\frac{S(t)}{t}\right).
\]
The singular expansions
$V(t)\overset{t\rightarrow  \alpha_2}{\sim}{1}/({1-\frac{t}{\alpha_2}})$,
$V'(t)\overset{t\rightarrow \alpha_2}{\sim}{1}/[{\alpha_2(1-\frac{t}{\alpha_2})^2}]$,
$S(t)\overset{t\rightarrow  \alpha_2}{\sim}{1}/({1-\frac{t}{\alpha_2}})$, and
$S'(t)\overset{t\rightarrow \alpha_2}{\sim}{1}/[{\alpha_2(1-\frac{t}{\alpha_2})^2}]$
yield the expansion
\[
G(t)\overset{t\rightarrow \alpha_2}{\sim}\frac{1}{\alpha_2(1-\frac{t}{\alpha_2})^2}.
\]
Consequently, at its dominant singularity $\alpha_2$, $U(z)$ satisfies
\begin{equation*}
U(z) = \frac{z}{(z-1)^2}\int_{0}^{z} G(t) \frac{(t-1)^2}{t} \, dt \overset{z\rightarrow \alpha_2}{\sim}\frac{\alpha_2}{(\alpha_2-1)^2}\int_{0}^{z}\frac{1}{\alpha_2(1-\frac{t}{\alpha_2})^2}\frac{(\alpha_2-1)^2}{\alpha_2} \, dt
\overset{z\rightarrow \alpha_2}{\sim}\frac{1}{1-\frac{z}{\alpha_2}}.
\end{equation*}
By applying eq.~(\ref{pillo1}) together with eq.~(\ref{x}), we finally have
\begin{equation}\label{etotc}
\mathbb{E}[T_n^2]=[z^n]U(z)\sim \alpha_2^{-n} \sim \mathbb{E}[T_nR_n] \sim \mathbb{E}[R_n^2] .
\end{equation}

From these calculations, we obtain the next result.
\begin{prop}\label{prop:varyule}
The variance of the total number $T_n$ of ancestral configurations in an ordered unlabeled history of size $n$ selected under the uniform distribution satisfies the asymptotic relation ${\mathbb V}[T_n] \sim {\mathbb E}[T_n^2] \sim {\mathbb V}[R_n] \sim \alpha_2^{-n}= (2.0449954971\ldots)^n$.
For increasing values of $n$, the correlation coefficient $\rho[T_n,R_n]$ between the total number $T_n$ of ancestral configurations and the number $R_n$ of root configurations converges to $1$, $\rho[T_n,R_n] \rightarrow 1$.
\end{prop}
\noindent \emph{Proof.}
For the variance, we observe that
${\mathbb V}[T_n] = \mathbb{E}[T_n^2] - \mathbb{E}[T_n]^2 \sim {\mathbb E}[T_n^2]$.
Indeed, from Proposition \ref{proppa}, we find $\mathbb{E}[T_n]^2 \bowtie (\alpha_1^{-2})^n$, where $\alpha_1^{-2} < \alpha_2^{-1}$, where $\alpha_2^{-1}$ is the exponential order of $\mathbb{E}[T_n^2]$, as in eq.~(\ref{etotc}). Also, ${\mathbb E}[T_n^2] \sim {\mathbb V}[R_n]$, because ${\mathbb E}[T_n^2] \sim {\mathbb E}[R_n^2] \sim {\mathbb V}[R_n]$ follows from eqs.~(\ref{etotc}) and (\ref{vardis-propyuv}).

Similarly, for the covariance between $T_n$ and $R_n$ we obtain
$${\rm Cov}[T_n,R_n] = \mathbb{E}[T_n R_n] - \mathbb{E}[T_n] \, \mathbb{E}[R_n] \sim {\mathbb E}[T_nR_n] \sim {\mathbb E}[T_n^2].$$
Indeed, from Proposition \ref{proppa} we have $\mathbb{E}[T_n] \, \mathbb{E}[R_n] \bowtie (\alpha_1^{-2} )^n$,
while from eq.~(\ref{etotc}), $\mathbb{E}[T_n R_n] \sim {\mathbb E}[T_n^2] \bowtie \alpha_2^{-n}$, with $\alpha_2^{-1} > \alpha_1^{-2}$.
Hence, the correlation coefficient between $T_n$ and $R_n$ is
\[
\rho[T_n,R_n] = \frac{{\rm Cov}[T_n,R_n]}{\sqrt{{\mathbb V}[T_n]} \sqrt{{\mathbb V}[R_n]}} \sim \frac{{\mathbb E}[T_n^2]}{{\mathbb E}[T_n^2]} = 1. \quad \Box
\]

\vskip .2cm
By the same argument of eqs.~(\ref{logN1}), (\ref{logN2}), and (\ref{eq:slutsky}), for uniformly distributed ordered unlabeled histories, the total number $T_n$ of ancestral configurations can be shown to follow an asymptotic lognormal distribution. In particular, the variables $({\log T_n-{\mathbb E}[\log T_n]})/{\sqrt{{\mathbb V}[\log T_n]}}$ and $({\log R_n-{\mathbb E}[\log R_n]})/{\sqrt{{\mathbb V}[\log R_n]}}$ for random ordered unlabeled histories of size $n$ converge asymptotically to standard normal distributions, where $\mathbb{E}[\log T_n] \sim \mathbb{E}[\log R_n] \sim 0.351 n$ and $\mathbb{V}[\log T_n] \sim \mathbb{V}[\log R_n] \sim 0.008 n$ (eq.~\ref{lognormalconst}).

\begin{figure}
\begin{center}
\includegraphics*[scale=0.40,trim=0 0 0 0]{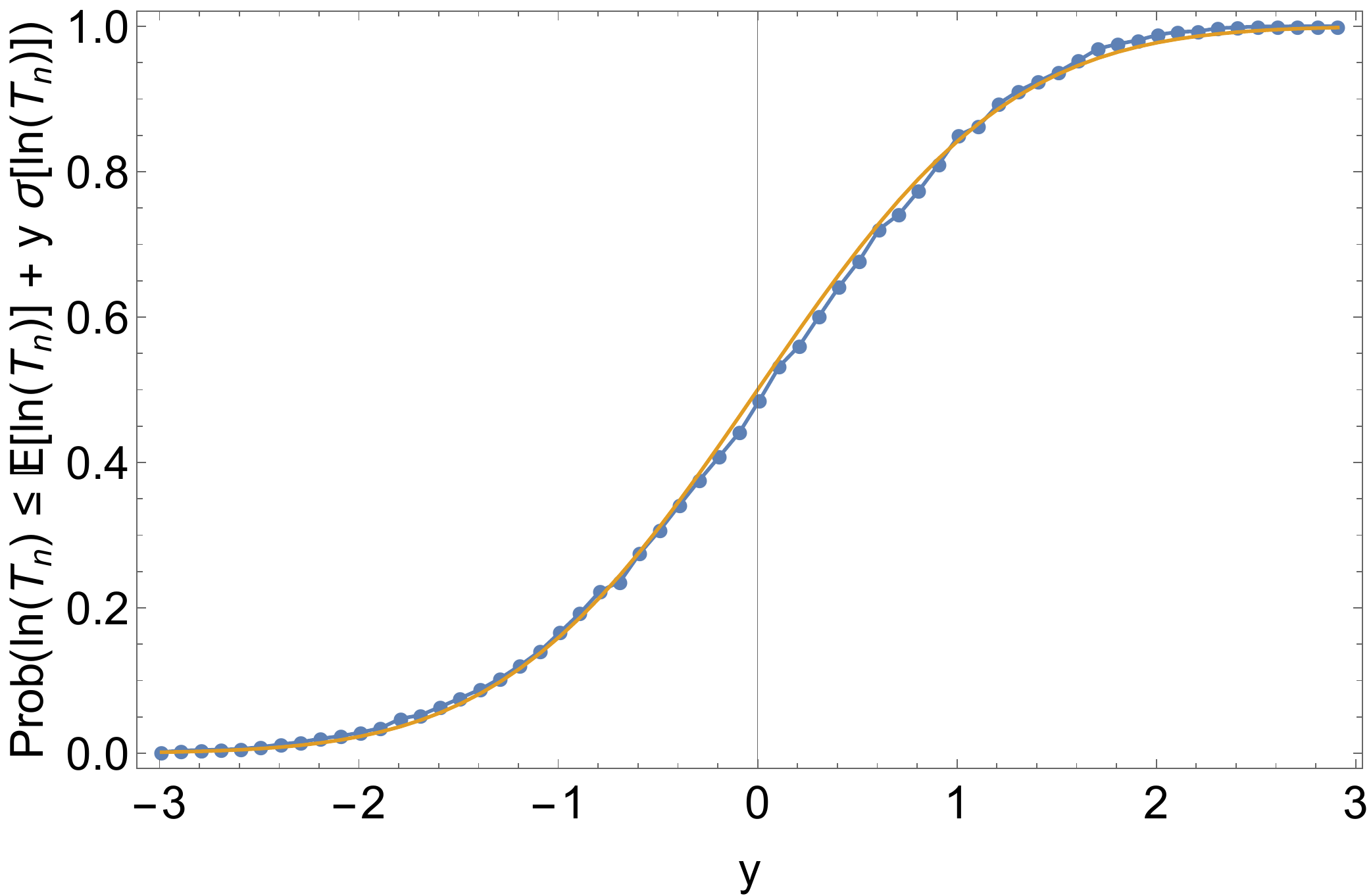}
\end{center}
\vspace{-.7cm}
\caption{{\small Cumulative distribution of the natural logarithm of the total number of configurations for uniformly distributed ordered unlabeled histories (or Yule-distributed labeled topologies) of size $n = 15$ (dotted line). For each $y \in [-3, 3]$ in steps of $0.1$, the quantity plotted is the probability that an ordered unlabeled history (or Yule-distributed labeled topology) with $n = 15$ chosen at random has a total number of configurations less than or equal to $\exp (\mathbb{E}[\log T_n] + y \sigma[\log T_n])$, where $\mathbb{E}[\log T_n]$ and $\sigma[\log T_n]=\sqrt{{\mathbb V}[\log T_n]}$ are respectively the mean and standard deviation of the logarithm of the total number of configurations for uniformly distributed ordered unlabeled histories (or Yule-distributed labeled topologies) with $n = 15$ taxa. The solid line is the cumulative distribution of a Gaussian random variable with mean 0 and variance 1.
}} \label{normalyule}
\end{figure}

In Fig.~\ref{normalyule}, we plot the cumulative distribution $\mathbb{P}[\log T_n \leq \mathbb{E}[\log T_n] + y \sqrt{{\mathbb V}[\log T_n]}]$ as a function of $y$, when ordered unlabeled histories of size $15$ are selected uniformly at random. To obtain the distribution, we can count total configurations for each unlabeled topology $t$ with 15 leaves, and then count the number (eq.~\ref{ouheq}) of ordered unlabeled histories having $t$ as tree shape. The figure illustrates the agreement between $({\log T_n-{\mathbb E}[\log T_n]})/{\sqrt{{\mathbb V}[\log T_n]}}$ and the standard normal distribution.

By the equivalence in Lemma \ref{lemcorr} between uniformly distributed ordered unlabeled histories and Yule-distributed labeled topologies, we summarize the results of this section.
 \begin{theo}\label{totttyu}
For a labeled topology of size $n$ selected at random under the Yule distribution, the mean and the variance of the total number $T_n$ of ancestral configurations grow asymptotically like
\begin{eqnarray}
\mathbb{E}[T_n] &\sim & \left(\frac{1}{1-e^{-2 \pi \sqrt{3}/9}}\right)^n,\\
\mathbb{V}[T_n] &\sim & ( 2.0449954971\ldots )^n.
\end{eqnarray}
Furthermore, the logarithm of the total number of ancestral configurations in a Yule-distributed labeled topology of size $n$ selected at random, rescaled as $({\log T_n-{\mathbb E}[\log T_n]})/{\sqrt{{\rm Var}[\log T_n]}} $, converges to a standard normal distribution, where ${\mathbb E}[\log T_n] \sim \mu n $ and ${\mathbb V}[\log T_n] \sim \sigma^2 n$, $(\mu, \sigma^2) \approx (0.351, 0.008)$.
\end{theo}

\section{Conclusions}

For a gene tree and species tree with matching labeled topology $t$ of size $n$ selected at random under the uniform and Yule probability models, we have studied the asymptotic distribution of the total number $T_n$ of ancestral configurations of $t$. By using techniques of analytic combinatorics, we have extended results of \cite{DisantoEtAl20} and \cite{DisantoAndRosenberg17}, where the focus was on the number $R_n$ of root ancestral configurations of $t$.


We have found that under both the uniform and Yule models, the total number of configurations has an asymptotic lognormal distribution, as was also demonstrated for the number of root configurations by \cite{DisantoEtAl20}. Furthermore, in Theorems \ref{totttuni} and \ref{totttyu}, we have shown that the mean and the variance of the total number of ancestral configurations grow like $\mathbb{E}[T_n] \sim \sqrt{6}(4/3)^n$ and $\mathbb{V}[T_n] \sim 5.050 \cdot 1.822^n$, for uniformly distributed labeled topologies, and like $\mathbb{E}[T_n] \sim 1.425^n$ and $\mathbb{V}[T_n] \sim 2.045^n$, when labeled topologies of size $n$ are selected under the Yule probability model. In particular, we observe that the mean total number of configurations is twice the mean number of root configurations under the uniform distribution for labeled topologies, with a correlation coefficient between $T_n$ and $R_n$ close to $0.9$ for $n$ large. For labeled topologies under the Yule distribution, the mean total number of configurations behaves asymptotically like the mean number of root configurations, with a correlation coefficient between $T_n$ and $R_n$ that approaches $1$ for increasing $n$. A summary appears in Table \ref{table:summary}.

That $T_n$ and $R_n$ have the same asymptotic growth under the Yule distribution on labeled topologies, and a correlation that approaches 1, is somewhat remarkable. $R_n$ tabulates ancestral configurations only at the root, whereas $T_n$ sums configurations across all $n-1$ internal nodes, \emph{including} the root. The correlation result indicates that under the Yule distribution, the configurations at non-root nodes contribute negligibly to the total.

\begin{figure}
\begin{center}
\includegraphics*[scale=0.40,trim=0 0 0 0]{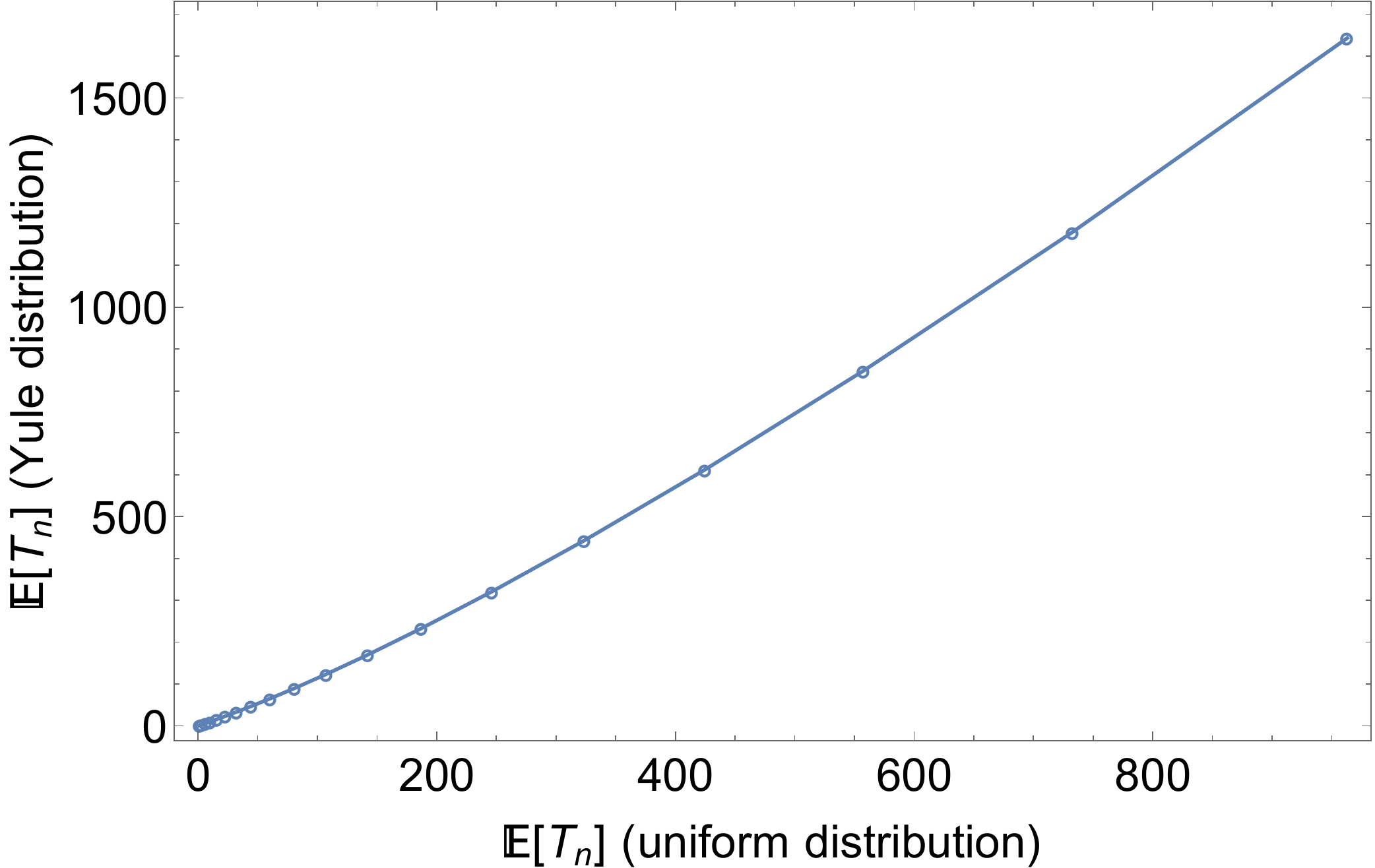}
\end{center}
\vspace{-.7cm}
\caption{{\small Mean total number of configurations of labeled topologies of size $n$ under the Yule and uniform distributions, for $2 \leq n \leq 20$. Values for the uniform distribution are computed by combining the recurrences in eqs.~(\ref{ern1}) and (\ref{etn1}); values for the Yule distribution are computed by combining the recurrences in eqs.~(\ref{recurr0}) and (\ref{recurr1}).}
} \label{yulevsuni}
\end{figure}

The difference in results for the uniform and Yule models suggests a correlation between tree balance and total configurations. Indeed, \cite{DisantoAndRosenberg17} suggested such a relationship for root configurations. A similar relationship might exist for total configurations; we find indeed that under the Yule model, which gives more weight to balanced topologies, the mean \emph{total} number of configurations grows faster than under the uniform model (Fig.~\ref{yulevsuni}).

\begin{figure}
\begin{center}
\includegraphics*[scale=0.80,trim=0 0 0 0]{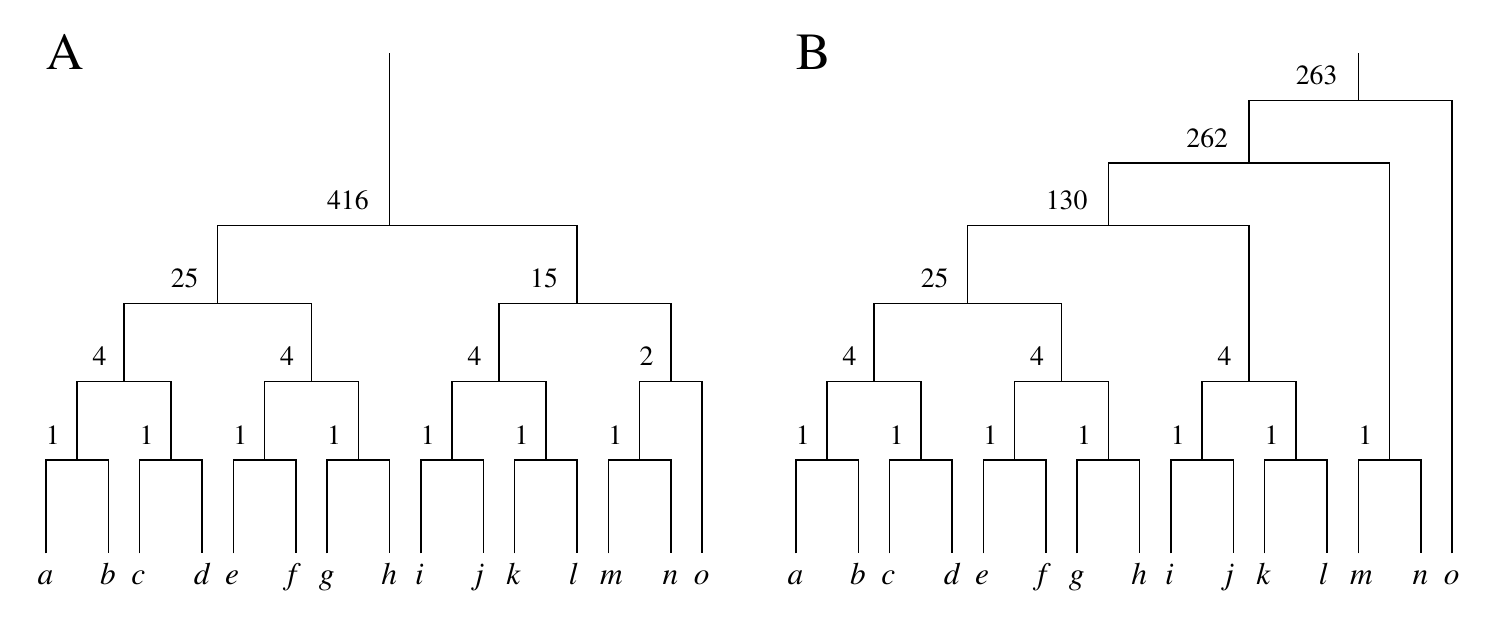}
\end{center}
\vspace{-0.85cm}
\caption{{\small The labeled topologies of size $n=15$ that (up to permutation of the labels) have maximal numbers of configurations. (A) Root  configurations (416). (B) Total configurations (690). Numbers of configurations appear at each internal node. The rightmost point in Fig.~\ref{totandroot} has coordinates $(\log 416, \log 477)$; the topmost point has coordinates $(\log 263, \log 690)$.}} \label{treestotconfig}
\end{figure}

Several directions naturally arise from our work. For instance, we did not characterize the labeled topologies of given size that have the largest total number of ancestral configurations. Section 4 of \cite{DisantoAndRosenberg17} described the recursive structure of labeled topologies with the maximal number of root ancestral configurations. However, as shown in Fig.~\ref{treestotconfig}, the number of root and total configurations do not necessarily attain their maximal values at the same labeled topology. We also did not consider non-matching gene trees and species trees. The non-matching case, in which the gene tree and species tree have different labeled topologies, merits further analysis, as a non-matching gene-tree labeled topology can have more total configurations than the topology that matches the species tree~\cite{DisantoAndRosenberg17}. It is of interest to see if techniques used in this article can be extended to derive distributional properties of the number of ancestral configurations when the gene tree and species tree differ in topology.

\begin{table}
\begin{center}
{\caption{\textcolor{black}{Summary of asymptotic equivalences for the number of root configurations and the total number of configurations under the uniform and Yule models on labeled topologies.}}\label{table:summary}}
\begin{tabular}{lllll}
\hline\hline
& \multicolumn{2}{c}{Uniform model} &
\multicolumn{2}{c}{Yule model} \\ \cline{2-5}
Quantity & Result & Reference & Result & Reference \\ \hline
$\mathbb{E}[R_n]$       & $\sqrt{\frac{3}{2}} \left(\frac{4}{3}\right)^n$ & eq.~(\ref{meandis-propyu}) & $\left(\frac{1}{ 1-e^{-2\pi\sqrt{3}/9} } \right)^n$ & eq.~(\ref{meandis-propyu}) \\
$\mathbb{E}[R_n^2]$     & $\sqrt{\frac{7(11 - \sqrt{2})}{34}}\left[\frac{4}{7 (8\sqrt{2} - 11) }\right]^n$ & eq.~(\ref{vardis-propyuv}) & $(2.0449954971\ldots)^n$ &  eq.~(\ref{vardis-propyuv}) \\
$\mathbb{V}[R_n]$       & $\sqrt{\frac{7(11 - \sqrt{2})}{34}}\left[\frac{4}{7 (8\sqrt{2} - 11) }\right]^n$ & eq.~(\ref{vardis-propyuv}) & $(2.0449954971\ldots)^n$ &  eq.~(\ref{vardis-propyuv}) \\
$\mathbb{E}[T_n]$       & $\sqrt{6} \left( \frac{4}{3} \right)^n$ & Prop.~\ref{meantotuni} & $\left(\frac{1}{ 1-e^{-2\pi\sqrt{3}/9} } \right)^n$ & Prop.~\ref{proppa} \\
$\mathbb{E}[T_n^2]$     & $\frac{2}{17}(15+11\sqrt{2})\sqrt{\frac{7(11-\sqrt{2})}{34}}\left[\frac{4}{7 (8\sqrt{2} - 11) }\right]^n$ & eq.~(\ref{ett2}) & $(2.0449954971\ldots)^n$ & Prop.~\ref{prop:varyule} \\
$\mathbb{V}[T_n]$       & $\frac{2}{17}(15+11\sqrt{2})\sqrt{\frac{7(11-\sqrt{2})}{34}}\left[\frac{4}{7 (8\sqrt{2} - 11) }\right]^n$ & Prop.~\ref{prop:varunif} & $(2.0449954971\ldots)^n$ & Prop.~\ref{prop:varyule}\\
$\mathbb{E}[T_n R_n]$   & $\left(1+\frac{\sqrt{2}}{2}\right)\sqrt{\frac{7(11-\sqrt{2})}{34}}\left[\frac{4}{7 (8\sqrt{2} - 11) }\right]^n$ & eq.~(\ref{ett1}) & $(2.0449954971\ldots)^n$ & Prop.~\ref{prop:varyule} \\
$\textrm{Cov}[T_n,R_n]$ & $\left(1+\frac{\sqrt{2}}{2}\right)\sqrt{\frac{7(11-\sqrt{2})}{34}}\left[\frac{4}{7 (8\sqrt{2} - 11) }\right]^n$ & Prop.~\ref{prop:varunif} & $(2.0449954971\ldots)^n$ & Prop.~\ref{prop:varyule} \\
$\rho[T_n,R_n]$         & $\frac{{1+\frac{\sqrt{2}}{2}}}{{\sqrt{\frac{2}{17}(15+11\sqrt{2})}}}$ & Prop.~\ref{prop:varunif} & 1 & Prop.~\ref{prop:varyule} \\
\hline\hline
\end{tabular}
\end{center}
We note the numerical values of recurring constants:
$1/(1-e^{-2\pi\sqrt{3}/9}) \approx 1.4253868277$,
$1-e^{-2\pi\sqrt{3}/9} \approx 0.70156394081$,
$4/[7 (8\sqrt{2} - 11)] \approx 1.8215272244$,
$7(8\sqrt{2} - 11)/4 \approx 0.5489898732$.
Constant $2.0449954971$ was evaluated in the Appendix of \cite{DisantoEtAl20}, and its reciprocal is $0.4889986317$.
\end{table}

\medskip
\noindent
{\footnotesize
{\bf Acknowledgments} This project developed from discussions at a workshop at the Banff International Research Station. Support was provided by National Institutes of Health grant R01 GM117590 to NAR, by grants MOST-107-2115-M-009-010-MY2 and MOST-109-2115-M-004-003-MY2 from Ministry of Science and Technology (MOST), Taiwan to MF, CYH and ARP, and by a Rita Levi-Montalcini grant from the Ministero dell'Istruzione, dell'Universit\`a e della Ricerca to FD.}




\end{document}